\documentstyle{article}
\newcommand{\1}{{{\mathchoice {\rm 1\mskip-4mu l} {\rm 1\mskip-4mu l}
{\rm 1\mskip-4.5mu l} {\rm 1\mskip-5mu l}}}}

\newcommand{\TWw}{{\Tilde{\Ww}}}

\newcommand{\tp}{{\tilde p}}
\newcommand{\ev}{{{\it ev}}}
\newcommand{\Cc}{{\cal C}}
\newcommand{\R}{{\bf R}}

\newcommand{\Q}{{\bf Q}}

\newcommand{\Z}{{\bf Z}}
\newcommand{\C}{{\bf C}}
\newcommand{\CP}{{\bf CP}}

\newcommand{\Aa}{{\cal A}}

\newcommand{\Mm}{{\cal M}}
\newcommand{\Nn}{{\cal N}}
\newcommand{\Hh}{{\cal H}}

\newcommand{\Ss}{{\cal S}}
\newcommand{\Ll}{{\cal L}}

\newcommand{\Uu}{{\cal U}}
\newcommand{\Ee}{{\cal E}}
\newcommand{\Xx}{{\cal X}}

\newcommand{\Ww}{{\cal W}}
\newcommand{\oMm}{{{\overline{\cal M}}\,\!}}
\newcommand{\pJ}{{\overline{\partial}_J}}

\newcommand{\ov}{\overline}
\newcommand{\oI}{{\ov{I}}}
\newcommand{\SO}{{\rm SO}}

\newcommand{\id}{{\rm id}}
\newcommand{\Tilde}{\widetilde}
\newcommand{\p}{{\partial}}
\newcommand{\al}{{\alpha}}

\newcommand{\be}{{\beta}}
\newcommand{\tbe}{{\tilde{\beta}}}
\newcommand{\tga}{{\tilde{\gamma}}}

\newcommand{\Om}{{\Omega}}
\newcommand{\om}{{\omega}}
\newcommand{\eps}{{\varepsilon}}
\newcommand{\de}{{\delta}}
\newcommand{\De}{{\Delta}}
\newcommand{\ga}{{\gamma}}
\newcommand{\Ga}{{\Gamma}}
\newcommand{\io}{{\iota}}
\newcommand{\ka}{{\kappa}}
\newcommand{\la}{{\lambda}}
\newcommand{\La}{{\Lambda}}
\newcommand{\si}{{\sigma}}

\newcommand{\tom}{\tilde\om}

\newcommand{\TS}{{\Tilde S}}
\newcommand{\ty}{{\tilde y}}
\newcommand{\tz}{{\tilde z}}
\newcommand{\tx}{{\tilde x}}

\newcommand{\TSs}{{{\Tilde \Ss}}}

\newcommand{\Symp}{{\rm Symp}}
\newcommand{\Diff}{{\rm Diff}}
\newcommand{\Ham}{{\rm Ham}}

\newcommand{\Aut}{{\rm Aut}}
\newcommand{\PD}{{\rm PD}}
\newcommand{\PSL}{{\rm PSL}}
\newcommand{\Si}{{\Sigma}}
\newcommand{\TU}{{\Tilde U}}
\newcommand{\TV}{{\Tilde V}}
\newcommand{\TLl}{{\Tilde \Ll}}
\newcommand{\TY}{{\Tilde Y}}
\newcommand{\TUu}{{\Tilde \Uu}}

\newcommand{\TZ}{{\Tilde Z}}
\newcommand{\TL}{{\Tilde L}}

\newcommand{\ttau}{{\tilde \tau}}
\newcommand{\th}{{\tilde h}}

\newcommand{\TE}{{\Tilde E}}
\newcommand{\TEe}{{\Tilde \Ee}}
\newcommand{\ts}{{\tilde s}}

\newcommand{\oMmn}{{\overline{{\cal M}}\,\!^\nu}}
\newcommand{\wMm}{{\widehat{\wMm}}}

\newcommand{\MS}{{\medskip}}
\newcommand{\SS}{{\smallskip}}

\newcommand{\NI}{{\noindent}}
\newcommand{\proof}[1]{\noindent{\bf Proof#1:\  }}
\newcommand{\QED}{\hfill$\Box$\medskip}
\newtheorem{theorem}{Theorem}[section]
\newtheorem{claim}[theorem]{Claim}
\newtheorem{thm}[theorem]{Theorem}
\newtheorem{cor}[theorem]{Corollary}
\newtheorem{lemma}[theorem]{Lemma}

\newtheorem{defn}[theorem]{Definition}
\newtheorem{prop}[theorem]{Proposition}
\newtheorem{rmk}[theorem]{Remark}
\newtheorem{remark}[theorem]{Remark}
\newtheorem{example}[theorem]{Example}

\begin{document}

\title{Quantum homology of fibrations over $S^2$}
\author{Dusa McDuff\thanks{Partially
supported by NSF grant DMS 9704825.} \\ State University of New York
at Stony Brook \\ (dusa@math.sunysb.edu)}

\date{May 12, 1999}
\maketitle

\vspace{1in}

\NI
{\bf Abstract}

This paper studies the (small) quantum homology and cohomology of 
 fibrations
$p:  P\to S^2$ whose structural group is the group of Hamiltonian
symplectomorphisms  of the fiber $(M,\om)$.  It gives a proof that
the rational cohomology splits additively as the vector space tensor product
$H^*(M)\otimes H^*(S^2)$, and investigates conditions under which 
the ring structure  also splits, thus generalizing work of
Lalonde--McDuff--Polterovich and Seidel. The main tool is a study of
certain operations in the quantum homology of the total space $P$ and 
of the fiber $M$, whose properties reflect the relations between
the Gromov--Witten invariants of  $P$ and  $M$.
In order to establish these properties we   further develop the language
introduced in [Mc3] to describe the virtual moduli cycle  (defined by
Liu--Tian, Fukaya--Ono, Li--Tian, Ruan and Siebert).

\vspace{1.2in} 

\NI
AMS classification number 53C15; 

\NI
key words: quantum cohomology, symplectic fibration,
Hamiltonian fibration, Gromov--Witten invariants

\newpage

\section{Main results}

\subsection{The homology of $P$}

This paper studies the quantum homology of Hamiltonian 
 fibrations
$$
p:  P\to S^2.
$$
A fibration $P\to B$ is called {\it symplectic} if its fiber is a symplectic manifold
$(M,\om)$ and its structural group  is the group of symplectomorphisms
$\Symp(M,\om)$.   Thus each fiber $M_b = p^{-1}(b)$  is equipped with a
well-defined symplectic form $\om_b$.  To say that $p$ is {\it Hamiltonian} 
means that  the structural group reduces further to the group of Hamiltonian
symplectomorphisms $\Ham(M,\om)$.    As noted by Seidel~[Sd1], when $B= S^2$ this
is equivalent to the existence of a symplectic form $\Om$ on $P$ that restricts
to $\om_b$ on each fiber.  Moreover, we can choose $\Om$ to be compatible with the
orientation of $P$ defined by the standard orientations on $S^2$ and on $(M,\om)$. 
Such forms $\Om$ are said to be {\it compatible} with $p$.  It is easy to see that the
set of all such forms
 is path connected.  This means in particular that the Gromov--Witten invariants of
$P$ do not depend on the choice of $\Om$.

We shall say that  
the rational cohomology $H^*(P) = H^*(P,\Q)$ of $P$ is {\em additively
split} if it is isomorphic (as  a vector space) to the tensor product
$H^*(M)\otimes H^*(B)$.  This happens if and only if the rational
homology $H_*(M)$ of $M$ injects into that of $P$.
The first aim of this paper is to give a proof of the following result that was
announced in Lalonde--McDuff--Polterovich~[LMP2].  It
has several interesting
corollaries related to the Flux conjecture for $(M,\om)$ that are fully
described there.  (Cf also [LMP1].)  The extent to 
which it generalises to Hamiltonian fibrations with other bases
 is discussed in the forthcoming paper~[LMP3].

\begin{thm}\label{main}  The rational cohomology of any Hamiltonian 
fibration with base $S^2$ is additively split.
\end{thm}

The reason for this is that each such  fibration gives rise to an operation on the
quantum homology $QH_*(M)$ of $(M,\om)$, whose existence  implies that the rational
homology $H_*(M)$ of $M$ injects into that of $P$.  This operation 
was first investigated by Seidel in [Sd1] and has the
form 
$$
a\mapsto {Q_\si} *_M a,\quad a\in H_*(M,\R),
$$
where $\si\in H_2(P,\Z)$ is a reference section, $ Q_\si\in QH_*(M)$
 represents the intersection with $M$ of all $J$-holomorphic sections of $P\to
M$, and $*_M$ is the  quantum product in $M$.  For more details see \S\S1.3,
 2.   Note
that we consider $Q_\si$ as an element of  $QH_*(M)$ rather than $H_*(M)$ in
order to keep track of the homology classes of the sections considered.  

The above theorem  was proved by Deligne when $p$ is a holomorphic  fibration
with K\"ahler total space.  
When the the loop  in $\Ham(M)$ that generates the
fibration comes from a circle action, it also follows from work of Kirwan on
equivariant cohomology.  (We show in Example~\ref{ex:loop}
 that not all loops have this form.)  More generally, Polterovich
pointed out that Blanchard [B] found an easy cohomological proof that applies
in the case of smooth fibrations with fiber  of ``hard Lefschetz type.''\footnote
{ This means that the map $\wedge [\om^k] : H^{n-k}(M,\R) \to H^{n+k}(M,\R)$
is an isomorphism for all $k$.}  Also,  it was proved in~[LMP2] in  cases
where the ``usual" theory of $J$-holomorphic spheres works, for example,
if $(M,\om)$ is spherically monotone.  In the present paper  the ``new" theory of
Gromov--Witten invariants (due to Fukaya--Ono~[FO], Li--Tian~[LiT], Liu-Tian~[LiuT],
Ruan~[R] and Siebert~[S])
 is used to extend the result to arbitrary $(M,\om)$.   In order to do this, 
we adapt 
this theory  to our fibered setting in order to show that the
operation we define has the requisite properties: see \S1.3, \S2, and \S4. \MS

Thus the additive structure of $H^*(P)$ is always very simple.
The story concerning its multiplicative structure is more complicated
since 
 one can consider both the standard cup product and also versions of
the quantum (or deformed) cup product.  For example, Seidel exploits
properties of the quantum product in his work on $\Symp({\CP}^m\times {\CP}^n)$
in~[Sd2].  He also pointed out\footnote {Private communication} that, if 
there is a fibered\footnote
{
 i.e. it restricts on each fiber to an $\om$-tame almost
complex structure: see Definition~\ref{def:Jfib} for a more precise
definition}
$\Om$-tame almost complex structure $J$ on $P$
 that admits  no $J$-holomorphic spheres in any fiber, then
$H^*(P)$ is isomorphic as a ring (under cup product) with the product of the
rings $H^*(S^2)$ and  $H^*(M)$.   For example, this would be the case if $\om
= 0$ on $\pi_2(M)$.   However, it is certainly not true  that the ring $H^*(P)$
always splits in this way: see Example~\ref{ex:PD}.  

We now state a generalization of this result with hypothesis framed in
terms of the vanishing of certain Gromov--Witten invariants.
 Given any symplectic manifold $(X,\om)$ we denote by
$n_X(v_1,\dots, v_k; B)$ the Gromov--Witten invariant that counts the number
of isolated $J$-holomorphic spheres in class $B\in H_2(X,\Z)$ that intersect
 representing cycles of the homology classes $v_1,\dots, v_k\in H_*(X)$.  (For
more details on these invariants, see \S\S1.3, 2, 4.)   
We shall denote the rational cohomology of $X$ by $H^*(X)$ and shall write
$\io:H_*(M)\to H_*(P)$ for the map induced by the inclusion of $M$ as a fiber.

\begin{thm}\label{ring} Let $(M,\om) \to P \to S^2$ be a 
fibration with structural
group $\Ham(M,\om)$, and suppose  that the Gromov--Witten invariants
$$
n_P(v_1,\dots, v_k; \io(B)),\quad v_j\in H_*(P),  B\ne 0, k \le 4,
$$
vanish.
Then  $H^*(P)$ is isomorphic as a ring under cup product with the
product of the rings $H^*(S^2)$ and $H^*(M)$. \end{thm}

 Other conditions under which the ring splits are stated in
Proposition~\ref{prop:ring2}.  However, the validity of the
suggestion made in~[Mc4] that the ring splits under the sole assumption that
the quantum product in $M$ is trivial is not yet clear.\MS

Our approach also allows us to find conditions under
which the nonsqueezing theorem holds
for  the  fibration  $P\to S^2$.
By this we mean the following.  Let $\Om$ be a symplectic 
form on $P$ that is compatible with the fibration $P\to S^2$
and define the {\it area} of $(P,\Om)$ to be the number $\al$ such
that $$
\frac{1}{(n+1)!}\int_P \Om^{n+1} \;=\; \frac {\al}{n!} \int_M
\om^{n}.
 $$
Thus, if the fibration $P\to S^2$ is symplectically trivial so that
$(P,\Om)$ is the product $(M\times S^2, \om \oplus\om_S)$, $\al$ is simply
 the area of the base $(S^2,\om_S)$.\footnote
{
Polterovich showed in [P] that the minimum of the  areas of all 
compatible symplectic forms on $P$ is an interesting measurement of the size of
the generating loop $\phi$ of $P$ (see \S 1.2) that is a kind of ``one-sided"
Hofer norm for $\phi$.}
  Then we will say
that {\it the  nonsqueezing theorem holds for the fibration} $p: (P, \Om) \to S^2$ if the
area $\al$ constrains the size of the balls that embed into $(P,\Om)$, i.e. if
$
\pi r^2 \le \al
$
whenever $ B^{2n+2}(r)$ embeds symplectically in $
(P,\Om)$.    By  considering the case when
$(P, \Om)$ is ${\CP}^2$ blown up at a point 
with one of the standard symplectic forms, it
is not hard to see that  the nonsqueezing
theorem does not hold for all $(P,\Om)$. 

The best result that our methods give about the nonsqueezing theorem  is
stated in \S3.4 as Proposition~\ref{prop:nonsq}.
Here is a corollary.

\begin{prop}\label{sque}  Suppose that the hypotheses of
Theorem~\ref{ring} hold.  Then, 
the nonsqueezing
theorem holds for the fibration $p:(P, \Om) \to S^2$.
\end{prop}

\subsection{Homomorphisms defined on $\pi_1(\Ham(M,\om))$}

In this section we discuss three
homomorphisms defined on
$\pi_1(\Ham(M,\om))$ that
 measure the cohomological twisting of the fibration
$P\to S^2$.  The most significant one is adapted from Seidel~[Sd1]
and is assembled from some of the ``horizontal" Gromov--Witten invariants of
$P$.

   First, recall
that any bundle over $S^2$ can be constructed by taking product bundles
over two discs and  gluing them by a clutching function $\phi = \{\phi_t\}
_{t\in [0,1]}$, viz:  $$
P_\phi = \left(D_+^2\times M\right) \cup_\phi \left(D_-^2\times M\right),
$$
where 
$$
\phi: (2\pi t,x)_+\mapsto (-2\pi t, \phi_t(x))_-.
$$
(The orientation on the base is induced from that on $D_+^2$.) Moreover  if we fix an
identification of one fiber of $P$ with $M$, we can normalise $\{\phi_t\}$ by requiring
that $\phi_0 = \id$ and get a loop that is well defined up to homotopy.
It is not hard to see that $P$ is a Hamiltonian fibration precisely when
$\phi$ is homotopic to a loop in  $\Ham(M,\om)$.      
Thus there is a bijective correspondence between loops $\phi \in
\pi_1(\Ham(M,\om))$ and  Hamiltonian  fibrations $P_\phi$ with one fiber
identified with $M$.

As noted in [LMP2], the manifold $P_\phi$ carries two canonical cohomology classes,
the first Chern class of the vertical
tangent bundle $$ 
c_\phi = c_1(TP_\phi^{vert})  \in H^2(P_{\phi},\Z),
$$
and the   coupling class  $u_\phi$, i.e.  the unique class
in $H^2(P_\phi,\R)$ such that 
$$
i^*(u_\phi) = [\om],\qquad u_\phi^{n+1} = 0.
$$
We will see in Lemma~\ref{le:glsect} that these classes behave well under
composition of loops.  

Since these classes are canonical, they give rise to
homomorphisms from $\pi_1(\Ham(M,\om))$.
The first, defined by Seidel in
 [Sd1] \S 10, has the form
$$
\oI_c\,:  \pi_1(\Ham(M,\om)) \to \Z/N\Z
$$
where $N$ is the minimal spherical Chern number, i.e. the smallest nonnegative
integer such that $c_1:H_2^S(M,\Z)\to\Z$ takes values in $N\Z$. 
(Here $H_2^S(M,\Z)$ be the spherical part of $H_2(M)$,
that is the image of the Hurewicz homomorphism $\pi_2(M)\to H_2(M,\Z)$.) In
our language,  $\oI_c$ is defined as  $$
\oI_c(\phi) = c_\phi(\si) \pmod {N},
$$
 where $\si$ is any section  class of $P_\phi$.
  This is well defined since $c_\phi(\si - \si')\in N\Z$ by
definition of $N$.   

The second homomorphism is determined by
 the Poincar\'e dual $\PD(u_\phi^n)\in H_2(P_\phi, \R)$. 
 Polterovich\footnote
{Private communication} pointed out that there is a canonical identification
between the quotients $H_2(P,\R)/H_2^S(P,\R)$ and 
$H_2(M,\R)/H_2^S(M,\R)$. Moreover, both groups are isomorphic to the second
homology group  $H_2(\pi_1(M),\R)$ of $\pi_1(M)$. (This follows by looking at the
homology spectral sequence of the fibration $\Tilde M \to M \to B(\pi_1(M))$ where
$\Tilde M$ is the universal cover of $M$ and $BG$ denotes the classifying
space of the group $G$.)  Hence there is a map   $$
\oI_u: \pi_1(\Ham(M))\to 
H_2(\pi_1(M),\R):\qquad \phi\mapsto  [\PD(u_\phi^n)]
$$
and it follows from Lemma~\ref{le:glsect} that this is a
homomorphism. By
Example~\ref{ex:PD} below, this homomorphism is not always trivial.

\begin{rmk}\rm
It would be interesting to  define  these homomorphisms
in a more intrinsic way that involved only $M$,
not $P_\phi$.  Since the definition of $\oI_c$ depends only on the
fact that the bundle $P\to S^2$ has a well defined fiberwise almost structure,
 $\oI_c$
extends to a homomorphism defined on the kernel of the
homomorphism $\pi_1(\Diff(M))\to \pi_1(\Aa)$ where $\Aa$ is the space of all
almost complex structures on $M$ that are homotopic to an $\om$-tame $J$.
Similarly, one can extend  the domain of  $\oI_u$ to the kernel of the
generalized flux homomorphism $F: \pi_1(\Diff(M)) \to H^1(M,\R)$
$$
F(\phi)(\ga) = \langle \om, \,{\phi_*(\ga)}\rangle,
$$
where $\phi_*(\ga)$ is the trace $T^2\to M: (s,t)\mapsto \phi_t(\ga(s))$
under $\phi$ of the loop $\ga(s)$.  This holds because $F(\phi) = 0$
precisely when there is a cohomology class in $H^2(P_\phi)$ that restricts to
$[\om]$, and this  is the  assumption underlying the definition of
the coupling class $u_\phi$ .
 A family $I_k$ of homomorphisms  closely related to $\oI_u$
was defined in Corollary 3.F of [LMP2] as follows:
$$
I_k: \pi_1(\Ham(M)\to \R:\qquad \phi\to \int_{P_\phi} (c_\phi)^k
(u_\phi)^{n+1-k}.
$$
In particular, $I_1(\phi) = c_\phi(\PD((u_\phi)^n).$  \QED
 \end{rmk}

Next, we define a version of the small quantum homology $QH_*(M, \La_R)$
with coefficients in a real Novikov ring $\La_R$. 
Set $c = c_1(TM)\in H^2(M,\Z)$.  Let $\La$ be the usual  Novikov ring of
the group $\Hh = H_2^S(M,\Z)/\!\!\sim$  with valuation $I_\om$ where $B\sim B'$
if $\om(B-B') = c(B-B') = 0$, and   let $\La_R$ be
the analogous (real) Novikov ring based on the group $\Hh_R =
H_2^S(M,\R)/\sim$.  Thus $\La$ is the completion of the group ring of $\Hh$
with elements of the form 
$$
\sum_{B\in \Hh} \la_B e^B
$$
where for each $\ka$ there are only finitely many nonzero
$\la_B\in \Q$ with $\om(B) > - \ka$. Similarly, the elements of
$\La_R$ are
$$
\sum_{B\in \Hh_R} \la_B e^B,
$$
where 
$\la_B\in \Q$ and there is a similar finiteness condition.  Set 
$$
QH_*(M)
= H_*(M)\otimes\La,\qquad QH_*(M,\La_R)
= H_*(M)\otimes\La_R.
$$
Then $QH_*(M)$ 
is $\Z$-graded with $\deg(a\otimes e^B) = \deg(a) + 2c_1(B)$.  It is best to
think of
$QH_*(M,\La_R)$ as $\Z/2\Z$-graded with 
$$
QH_{\ev} =  
H_{\ev}(M)\otimes\La_R, \quad QH_{odd} =  
H_{odd}(M)\otimes\La_R.
$$

Recall that the quantum intersection product 
$$
a*_Mb\in QH_{i+j - 2n}(M), \qquad \mbox {for }\; a\in H_i(M), b\in H_j(M)
$$
 is defined as follows:
$$
a*_Mb = \sum_{B\in \Hh} (a*_Mb)_B\otimes e^{-B},
$$ 
where  $(a*_Mb)_B\in H_{i+j- 2n+2c(B)}(M)$ is defined by the requirement that
\begin{equation}\label{eq:qm}
(a*_Mb)_B\,\cdot_M\, c = n_M(a,b,c;B) \quad\mbox{ for all }\;c\in H_*(M).
\end{equation}
Here we have written $\cdot_M$ for the usual  intersection
pairing on $H_*(M)$.  Thus
$a\cdot_M b = 0$ unless $\dim(a) + \dim (b) = 2n$ in which case it is the
algebraic number of intersection points of the cycles. The product  $*_M$ is
extended to $QH_*(M)$ by linearity over $\La$, and is associative. Moreover,
it  preserves the grading if we set $$ {\rm deg\,}
 c\otimes e^B = \dim (c) + 2c_1(B).
$$
(The above formula is consistent with [Sd1]  and  is appropriate for quantum
homology rather then cohomology.)  The product clearly extends also to $\La_R$. Note
here that when defining $a*_Mb$ we still sum over classes $B\in \Hh$ (and not $B\in
\Hh_R$), since $J$-holomorphic spheres can only represent integral classes.

 With respect to this product $*_M$, both versions of quantum homology are
graded commutative rings with unit $\1 = [M]$. Further, the invertible
elements in $QH_{\ev}(M,\La_R)$ form a group  $QH_{\ev}(M,\La_R)^\times$
that acts on   $QH_*(M, \La_R)$ by  quantum multiplication. 

The following result is proved in \S3. It is a mild generalization of a result of
Seidel~[Sd1].  

\begin{theorem}\label{repr}   For each symplectic manifold
$(M,\om)$ there is a natural
homomorphism 
$$
\rho:  \pi_1(\Ham(M,\om)) \to QH_{\ev}(M,\La_R)^\times.
$$
 \end{theorem}

In fact  $\rho(\phi) = Q_\si$, where $Q_\si$ is made from the
$J$-holomorphic sections of $P$ as described just after the statement of
Theorem~\ref{main} and $\si = \si_\phi$ is a suitable reference section
that is defined in Lemma~\ref{le:si}.  The homomorphism $\rho$
 is determined by  Gromov--Witten invariants of the form
$$
n_P([M], [M], \io(a);\si)
$$
where $\io: H_*(M)\to H_*(P)$ is the inclusion and $\si$ is a section class.

Seidel showed that $\rho$ is often nontrivial.  For example, 
$\pi_1(\Ham(S^2,\om)) \cong \pi_1(\SO(3)) = \Z/2\Z$ is generated by a loop $\phi$
that consists of a full turn about some axis and one can check that 
$$
\rho(\phi) = [pt]\otimes e^{A/2}, \quad\mbox{where }\; A = [S^2].
$$
Observe also that if $[\om]$
is rational one can replace $\Hh_R$ by the corresponding rational group
$\Hh_Q$. 

We will see in \S2.1 that the very existence of $\rho$ is enough to prove
that $H^*(P_\phi)$ is additively split for every $\phi$.  However, it is not clear
if the condition  $\rho(\phi) = \1$  guarantees that
$H^*(P_\phi)$ splits as a ring, though so far there are no counterexamples
either.  Statement (i) in Proposition~\ref{prop:ring2} is as far as we have been
able to go in this direction.

\subsection{The fiberwise quantum homology of $P$}

Most of the proofs of the above results are based on the
relations between the Gromov--Witten invariants of $M$ and those of $P_\phi$.
Here we will state these relations explicitly. 
In order to put them in a nice form we will suppose in this section
that Theorem~\ref{main} is proven.  This implies that the Wang sequence for the
homology of $P_\phi$ degenerates into the short exact sequence:
$$
0\to H_{*+2}(M, \R)\stackrel {\io}{\to} H_{*+2}(P_\phi, \R) 
\stackrel{\cap[M]}\to
H_{*}(M, \R)\to 0,
$$
where  $\io$ is the map induced on
homology by the inclusion. 
Let $s: H_{*}(M, \R)\to  H_{*+2}(P_\phi, \R)$ be any splitting
of this sequence. Note the identity:
$$
s(a)\cdot_P \io(b) = a\cdot_M b.
$$
Besides the $3$ point Gromov--Witten invariant mentioned above, we consider the
$4$ point invariant in which the cross ratio of the marked points on the domain is
fixed.  We denote this by
$$
n_{P,\chi}(v_1,v_2,v_3,v_4; D).
$$
Intuitively it counts the number of $J$-holomorphic maps $h:S^2\to P$ in class $D$ 
such that $h(z_i)\in \nu_i, $ $ i = 1,\dots,4$, where $\nu_i$ is a cycle
representing the class $v_i$ and where the cross ratio of the $4$ points
$z_1,\dots,z_4$ on $S^2$ is fixed.  The following
 important identity is equivalent to the associativity of quantum
multiplication:
\begin{eqnarray}\label{eq:assoc4}
n_{P,\chi}(v_1,\dots, v_4; A) = \sum_{A= A_1+A_2, \,\al} 
 n_{P}(v_1, v_2, e_\al'; A_1)\,\cdot\,
n_{P}(f_\al', v_3, v_4; A_2)
\end{eqnarray}
where $\{e_\al'\}_\al$ is a basis for the homology $H_*(P,\Q)$, and
$\{f_\be'\}_\be$ is the dual basis with respect to the intersection pairing. 
 Observe also that often the class $B$ or
$D$ is not an individual homology class, but rather is an equivalence class of
such: see for example equation~(\ref{eq:qm}) above.  Whenever it is not clear
from the context
 what is
meant, either interpretation will do.

The following proposition expresses the
compatibility between the invariants of $P_\phi$ and $M$.  A similar statement 
(applied to the fibered space $\Xx$ with fibers $X_t$) forms a  crucial step in the
proof of Theorem~\ref{main}: see \S~\ref{ssec:comp}.  
A spherical class $\si\in H_2^S(P; \Z)$ is called a section class if $\si\cdot [M] =
1$.

\begin{prop} \label{prop:GW}  For any classes $a,b,c\in H_*(M), v,w\in H_*(P),
B\in H_2(M,\Z)$,  any splitting $s$ and any section class $\si$ we have:\MS

\NI{\rm (i)} $ n_P (\io(a),\io(b), v; \io(B)) = 0$.
\SS

\NI{\rm (ii)}  $n_P(\io(a), v,w;\io(B)) = n_M(a, v\cap [M], w\cap [M]; B)$.  In
particular $$
\begin{array}{lcl}
n_P([M], v, w;\io(B)) & = & 0 \;\mbox{ if }\; B\ne 0;\\
n_P(\io(a),s(b), s(c); \io(B))  & = & n_M(a,b,c\,; B).
\end{array}
$$
\SS

\NI
{\rm (iii)}\,$
n_P(\io(a), \io(b), v\cap [M]; \si)  = n_{P,\chi}(\io(a), \io(b), v, [M]; \si).
$

\NI
{\rm (iv)}\,\,\,\,\,\, $ n_P(v, \io(a), \io(b); \si)$
\begin{eqnarray*}
& =& \sum_{B,i} n_P(v, [M],
\io(e_i); \si -\io(B))\cdot n_M(f_i, a,b; B)\;  \\
& &  \qquad\quad + n_P([M], [M], \io(e_i); \si - \io(B)\cdot
n_M(f_i,v\cap [M],a,b;B)
\end{eqnarray*}
where $\{e_i\}$ is a basis for $H_*(M)$ with dual basis $\{f_i\}$.
\end{prop}

Although this proposition is very close to well known properties of Gromov--Witten
invariants,  there
do not seem to be any relevant references in the literature.  
Any expert in the
field would no doubt be able to supply a proof: we give one  in
\S 4.4.2.  (iii) and (iv) are especially noteworthy.  The latter is an analog of the
associativity rule~(\ref{eq:assoc4}) and is derived by moving the representing
cycles for $\io(a), \io(b)$ into the same fiber instead of  moving two of the 
marked points together.  In this case we count stable maps that have two
components, one a section in $P$ and the other a vertical curve (i.e.  lying in a
fiber.) The first sum corresponds to the case when the section meets $v$ and
the vertical component meets $a,b$, and the second to the case when the
vertical component meets all three cycles $a,b,v\cap [M]$.\MS

Next we introduce a fiberwise quantum product in $P$. 

\begin{defn}\rm Let $QH_*^V(P,\La) = H_*(P)\otimes\La$, where $\La = \La(M)$
is as before.  Given $u,v\in H_*(P)$, we define the vertical quantum product
in $P$ by setting  $u*_V v = \sum_B(u*_V v)_B\otimes e^{-B}$ where 
$$
(u*_V v)_B\cdot_P w = n_P(u,v,w;\io(B)).
$$
 \end{defn}

The next result exhibits the relations between the products $*_M$ and $*_V$.

\begin{lemma}\label{le:v} For all $a,b \in H_*(M)$ and any splitting $s$, 
we have\MS

\NI
{\rm (i)} $\io(a) *_V \io(b) = 0$,\SS

\NI
{\rm (ii)} $s(a)*_V \io(b) = \io(a*_M b)$,\SS

\NI
{\rm (iii)}  $(s(a) *_V s(b)) \cap [M] = a*_M b$.
\end{lemma}
\proof{} This can be proved by using
Proposition~\ref{prop:GW} and  the fact that the class $\io(a)$ is determined
by the  identities
$$
\io(a)\cdot_P s(b) = a\cdot_M b,\quad \io(a)\cdot_P \io(b) = 0,\qquad b\in
H_*(M).
$$
Details are left to the reader. \QED

The following  quantum analog of the Wang sequence in homology
is an immediate consequence.

\begin{prop} \label{prop:wang} The sequence 
$$
0\to QH_{*+2}(M, \La_R)\stackrel {\io}{\to} QH^V_{*+2}(P_\phi, \La_R) 
\stackrel{\cap[M]}\to
QH_{*}(M, \La_R)\to 0.
$$
is a short exact sequence of $QH_{*}(M, \La_R)$-modules, where the element
$a\in QH_{*}(M, \La_R)$  acts on  $QH_{*}(M, \La_R)$  via $b\mapsto
a*_M b$ and on $QH^V_{*+2}(P_\phi, \La_R)$  via 
$$
v\mapsto   s(a) *_V v.
$$
\end{prop}

Note that the module structure on the middle term depends, in
principle, on the choice of splitting $s$.  
There is a dual version of this for cohomology stated
in \S3.2 below that generalizes a result in 
Seidel~[Sd2]. In this paper we have chosen to work mostly with homology
since this formulation seems to display the geometry most transparently.
However, as is evidenced by Seidel's work, in some ways the structure on
cohomology is easier to understand. \MS

As mentioned above, the proof of Theorem~\ref{main} is based on the properties of an
 operation 
$$
\Psi_{\phi,\si}: QH_*(M)\to QH_*(M)
$$
 that is associated to the fibration $P_\phi$.  (Here $\si$ is a reference section class in
$H_2(P_\phi,\Z)$.)   One can define  this operation  in terms of the 
(first order) {\em horizontal
part} $*_{H,\si}$ of the full quantum  multiplication $*_P$ on $P_\phi$
as follows. 

\begin{defn}\rm  Given a reference section $\si$, identify the set of all
section classes in $H_2(P,\Z)$ with $H_2(M,\Z)$ by writing them as $\si + \io(B)$. 
Then,
 define  the multiplication $*_{H,\si}$  on $QH_*^V(P)$ by setting
  $u*_{H,\si} v = \sum_B(u*_{H,\si} v)_B\otimes e^{-B}$ where 
$$
(u*_{H,\si} v)_B\cdot_P w = n_P(u,v,w;\si + \io(B)).
$$
 \end{defn}

Thus $*_{H,\si}$  is exactly that 
part of $*_P$ coming from the section classes,
just as $*_V$ is the part coming from the fiber classes.
One way to define $\Psi_{\phi,\si}$ is as:
$$
\Psi_{\phi,\si}(a) = \left([M]*_{H,\si} \io(a)\right)\cap [M].
$$
We will see in \S3.1 that the homomorphism $\rho$ of
Theorem~\ref{repr} can be defined by setting
$$
\rho(\phi) = ([M]*_{H,\si_\phi}[M])\cap [M]
$$
for a suitable choice of (generalized)  section $\si_\phi$. Moreover, 
in the situation considered in
Theorem~\ref{ring} there is a section class $\si_A$ with $c_1(\si_A) = 2$ such
that the map 
$$
a\mapsto s_A(a) = [M] *_{H,\si_A} a
$$
splits the homology ring.

\MS

It would of course be very interesting to understand the full quantum
(co)homology of $P_\phi$.  In particular, to what extent is it determined by
$QH^*(M)$ and  Seidel's homomorphism $\rho$?  Givental gives in [Gi] some
hints as to what an answer might be, at least in the case when $\phi$ comes
from a circle action.  The work in \S3.3 can be considered as the beginning of
an attempt to tackle this question: in particular in
part (ii) of Proposition~\ref{prop:ring2} 
we show that assumptions on the Gromov--Witten invariants of $M$ and on
$\rho$ do give information on the (undeformed) cup product on $H^*(P)$.
The easiest case is when $P$ is the product  $M\times S^2$. Kontsevich and
Manin~[KM] showed (among many other things) that 
$$
QH^*(M\times S^2) = QH^*(M)\otimes QH^*(S^2)$$
when $M$ is a projective algebraic manifold.   The next result extends this
to all symplectic $M$.  This is surely well known, and we include a sketch of
the proof   for the sake of completeness:  see
Proposition~\ref{prop:prodt}.

\begin{prop}\label{pr:prod}  For all closed symplectic manifolds $M$, 
$$QH^*(M\times S^2) = QH^*(M)\otimes QH^*(S^2).
$$
\end{prop}

\MS

\NI
{\bf Final remarks.} \,\, This paper has
been greatly influenced by work of Seidel [Sd1] and Piunikhin--Salamon--Schwartz
[PSS]. 
   As in~[PSS], it is possible to define operations analogous to
$\Psi_{\phi,\si}$  on  fibrations over
other compact Riemann surfaces besides the sphere.    Our approach 
to Gromov--Witten invariants could
also be adapted to deal with general symplectic fibrations $X\to B$ with
arbitrary (non symplectic) base manifold $B$. 
This would lead to generalizations of the fiberwise (or family) Gromov--Witten
invariants that were considered in the semi-positive case by 
Le--Ono in~[LO], and would 
allow one to apply Seidel's ideas in~[Sd2] to any symplectic manifold: cf.
Remark~\ref{rmk:sei}.

 \MS

\NI
{\bf Acknowledgements}  
  The author has profited from many conversations about
this details of this paper with Polterovich and Lalonde.  Any remaining mistakes are
her own responsibility.  She also wishes to thank Seidel for stimulating
conversations, Salamon for suggesting that the virtual moduli cycle is a
branched manifold, Siebert, Fukaya and Tian for useful
discussions about general Gromov--Witten invariants, and ETH, Zurich, for its
hospitality during part of the work on this paper.

\section{Proof of Theorem~1.1}

In this section we define the operation $\Psi_{\phi,\si}$ on quantum homology
  and discuss its main properties.  Since these are based on facts
about  Gromov--Witten invariants, \S2.2 is devoted to a preliminary discussion of
these.  In \S2.3 we prove  the basic propositions~\ref{prop:const}
and~\ref{prop:comp} modulo Claim~\ref{cl} and Lemma~\ref{le:xx}.

\subsection{Quantum operations}

We begin by defining the operation on quantum homology associated
to $P_\phi$.   We will start from a definition 
that looks a little different from that given in
 \S1.3 since our formulations there  assumed the validity of
Theorem~\ref{main}.
 Two sections $\si, \si'$ of $P_\phi$ are said to be {\it equivalent} if
$$
u_\phi(\si) = u_\phi(\si'),\quad c_\phi(\si) = c_\phi(\si').
$$

\begin{defn}\label{def:psi}\rm 
Given a loop of Hamiltonian diffeomorphisms $\phi$ on $M$,
and an equivalence class of sections $\si$ of $P_{\phi}$
with $d = 2c_\phi(\si)$, we define the $\La$-linear map
$$
\Psi_{\phi,\si}: \;\; QH_*(M) \to QH_{*+d}(M)
$$
as follows.   For $a \in
H_*(M)$, $\Psi_{\phi,\si}(a)$ is  the class in $QH_{*+d}(M)$ whose
intersection with $b\in H_*(M)$ is given by: 
$$
  \Psi_{\phi,\si}(a) \cdot_M
b = \sum_{B\in \cal H} n_P(\io(a),\io(b); \si + \io(B))\otimes e^{-B}.
$$
\end{defn}

Here $n_P(v,w;\si')$ is the Gromov--Witten invariant
which counts isolated $J$-holomorphic stable curves in $P_\phi$ of genus $0$
and with two marked points, that
 represent the equivalence class $\si'$ and whose marked points go
through  given generic
representatives of the classes $v$ and $w$ in $H_*(P_{\phi})$.\footnote
{
Observe that $n_P(v,w;\si')$  is identical to the
$3$-point invariant $n_P(v,w,[M];\si')$ and we will use these notations
interchangeably. See equation~(\ref{eq:23}) in \S2.2.2.}
 Note that,
by Gromov compactness, there are for each
given energy level $\ka$  only finitely many section classes $\si'=\si +\io(B)$
with $\om(B) \le \ka$ that are represented by $J$-holomorphic curves
in $P_\phi$.  Thus
 $\Psi_{\phi,\si}(a)$ satisfies the finiteness condition for
elements  of $QH_*(M,\La)$.

For reasons of dimension, $n_P(v,w; \si') = 0$ unless 
\begin{eqnarray}\label{eq:dim1}
2c_\phi(\si') + \dim (v) + \dim (w) = 2n.
\end{eqnarray}
(Recall from \S1.2 that $c_\phi$ is the vertical Chern class, not the full first
Chern class $c_1$ of $P$.  Thus $c_\phi(\si') = c_1(\si') - 2$.) Further
$$
\Psi_{\phi,\si}(a) = \sum a_{\si,B}\otimes e^{-B},\quad a_{\si,B}\in H_*(M),
$$
where $a_{\si,B}\cdot_M b = n_P(\io(a), \io(b); \si + \io(B))$, and
\begin{eqnarray}\label{eq:dim}
\dim (a_{\si,B}) = \dim(a) + 2c_\phi(\si+\io(B))  = \dim (a) + 2c_\phi(\si) + 2c(B).
\end{eqnarray}
 Observe also that
$$
\Psi_{\phi,\si + A} = \Psi_{\phi,\si}\otimes e^{A}.
$$

The rest of this section is concerned with  two basic propositions.
Here is the first.

\begin{prop}\label{prop:const} If $\phi$ is the constant loop $*$ and $\si_0$
 is the class of the flat
section $[pt\times S^2]$ in $P_* = M\times S^2$, then $\Psi_{*,\si_0}$ is the
identity map.
\end{prop}

The second deals with the behavior of $\Psi_{\si,\phi}$ under composition of loops,
and before stating it we need some preparation.
One can represent the  composite $\psi*\phi$ 
 either by the product $\{\psi_t\circ\phi_t\}$ or by the concatenation
of loops.

\begin{lemma}\label{le:comp}  The bundle $P_{\psi*\phi}$ is diffeomorphic to
 the fiber sum $P_\psi\#_M P_\phi$.
\end{lemma} 
\proof{}  Recall that
$$
P_\phi = \left(D_+^2\times M\right) \cup_\phi \left(D_-^2\times M\right),
$$
where 
$$
\phi_t: (2\pi t,x)_+\mapsto (-2\pi t, \phi_t(x))_-.
$$
Let $M_{\phi,\infty}$ denote the fiber at $0\in D_2^-$ in $P_\phi$ and
$M_{\psi,0}$ the fiber at $0\in D_2^+$ in $P_\psi$. Cut out open product
neighborhoods of each of  these fibers and then glue the
complements by an orientation reversing symplectomorphism of the
boundary.   The resulting space may be realised as
$$
\left(D_2^+\times M\right)\;\cup_{\al_{\phi, -1}}\;\left(S^1\times
[-1,1]\times M\right)\;\cup_{\al_{\psi,1}}  \;\left(D_2^-\times M\right),
 $$
where
$$
{\al_{\phi, -1}}(2\pi t,x) = (2\pi t, -1,\phi_t(x)),\quad 
{\al_{\psi, 1}}(2\pi t, 1,x) = (-2\pi t,\psi_t(x)),
$$
and this may clearly 
be identified with  $P_{\psi*\phi}$. \QED

Using the above notation, let 
$$
\begin{array}{lcl}
V_\phi & = &\left(D_2^+\times M\right)\;\cup_{\al_{\phi, -1}}\;\left(S^1\times
[-1,1/2)\times M\right),\\\vspace{.1in}

 V_\psi & = &
\left(S^1\times
(-1/2,1]\times M\right)\;\cup_{\al_{\psi,1}}  \;\left(D_2^-\times M\right)
\end{array}
$$
The following lemma was proved in [LMP2].

\begin{lemma}\label{le:glsect} The classes $u_{\psi*\phi}$ and $c_{\psi*\phi}$ are
compatible with the decomposition $
P_{\psi*\phi}= V_\psi\cup V_\phi$ in the sense that their
 restrictions to $V_\psi\cap V_\phi = (-1/2,1/2)\times S^1 \times M$ equal
the pullbacks of $[\om]$ and $c_1(TM)$ under the obvious projection
$V_\psi\cap V_\phi\to M$.
\end{lemma}

Given  sections $\si$ of $P_\phi$ and  $\si'$ of
$P_\psi$ that agree over a small neighborhood of the fibers 
$M_{\phi,\infty}, M_{\psi,0}$, we write  $\si'\#\si$ for
the union of these sections in the fiber sum $P_\psi\#P_\phi = P_{\psi*\phi}$.
Here is the second basic result.

\begin{prop}\label{prop:comp}  For any sections   $\si$ of $P_\phi$ and  $\si'$
of $P_\psi$
$$
\Psi_{\psi,\si'}\circ \Psi_{\phi,\si} = \Psi_{\psi*\phi,
\si'\#\si}.
$$
\end{prop}

\begin{cor}\label{cor:comp}  $\Psi_{\phi,\si}$ is an isomorphism for all loops
$\phi$ and sections $\si$. 
\end{cor}
\proof{}   Observe that there is a fiberwise diffeomorphism from $P_\phi$ to
$P_{\phi^{-1}}$  that is a symplectomorphism on the fibers and covers an 
orientation reversing map on the bases.  Given a section $\si$ of $P_\phi$, let $-\si$
denote the image of this section in $P_{\phi^{-1}}$.  Then it is easy to check that
$\si\#(-\si)$ is isotopic to the flat section in the trivial bundle  $P_\phi\#P_{\phi^{-1}}
= P_*$. Hence, by Propositions~\ref{prop:const} and~\ref{prop:comp}, the
inverse   of   $\Psi_{\phi,\si}$ is $\Psi_{\phi^{-1}, -\si}.$\QED

\NI
{\bf Proof of Theorem~\ref{main}:}

The Wang sequence for the fibration $p: P_\phi\to S^2$
$$
\dots \to H_*(M) \stackrel{\io}{\to} H_*(P_\phi)\stackrel{\cap [M]}{\to}
H_{*-2}(M) \to H_{*-1}(M) \stackrel{\io}{\to}\dots
$$
implies that the spectral sequence for  $H_*(P_\phi)$ splits if and only if 
the map $\io:H_*(M)\to H_*(P_\phi)$ is injective.  Now Gromov-Witten
invariants are linear in each variable. Thus if $\io(a)=0 $ for some $a \neq 0$, then
$\Psi_{\phi,\si}(a) = 0$, a contradiction with the fact that $\Psi_{\phi,\si}$ is an 
isomorphism. \QED

One can almost prove  Propositions~\ref{prop:const} 
and~\ref{prop:comp}  using the axioms of general Gromov--Witten invariants. 
The missing step is to show that  these invariants can be  calculated
in a set-up that is adapted to our fibered setting.   In the rest of this section, we outline
proofs of these propositions  so that readers can see what is 
involved. Initially, we will use
Siebert's approach to Gromov--Witten invariants since in some ways this is
conceptually the simplest. 
 
\begin{example}\label{ex:loop}\rm  Since Theorem~\ref{main} is already
known in the case when the loop $\phi$ comes from a circle action, 
we give an example to show that not all loops in $\pi_1(\Ham(M))$ arise this
way.   Take $(M,\om) = (S^2\times S^2, (1+\la)\tau_1\oplus
\tau_2)$, where $0 < \la < 1$ and $\tau_i$ is an area form on the $i$th sphere
with total area $1$.   Put
$G^\la = \Ham(M,\om)$.  Given loops $u_t, v_t, t\in S^1,$ in
$G^\la$, let $u\times v$ denote the map $$
u\times v: T^2\to G^\la:\quad (s,t)\mapsto u_s v_t.
$$
First observe that, if $h_t, t\in S^1$, is homotopic
to a circle action then $h\times h$  is nullhomotopic.  Indeed, we can assume
that $h_sh_t = h_{s+t}$ so that the map $h\times h$ extends to $S^1\times D^2$
where $\p D_2$ is attached to the circles $s+t = const$.  Now think of
$S^2\times S^2$ as a fibration over the first sphere and let $\si_S, \si_N$ be
two sections of the form $S^2\times \{pt\}$, $\si_-$ be the anti-diagonal and
$\si_+$ the diagonal.  Further, for $t\in S^1$, let $u_t,$ (resp. $v_t$,) be a full
rotation in the fibers that fixes the sections $\si_N,\si_S$ (resp. $\si_\pm$.)
Because $\la> 0$ we can represent  $u_t, v_t$ by loops in $G^\la$. 
Moreover, it follows from  Anjos's calculation of the Pontriagin
ring $H_*(G^\la,\Z/2\Z)$ in [An] that 
$$
[u\times v] \ne [v\times u] \,\, \in\,\, H_2(G^\la, \Z/2\Z).
$$
Thus $(u+v)\times (u+v) = u\times v + v\times u\ne 0$, so that the element
$u+v$ is not homotopic to a circle action.
 \end{example}

\subsection{Gromov--Witten invariants}

For the convenience of the reader we begin by recalling the definition of stable maps.
For further information, see any of the papers~[FO], [LiT], [LiuT], [R] or [S].  Useful
background may be found in~[MS1] or [RT].

  \subsubsection{Stable  maps and curves}

Given a symplectic manifold $(P,\om)$ with compatible almost complex
structure $J$, denote by $\Mm_{0,k}(P,J,A)$  the space of
(unparametrized) $J$-holomorphic spheres  in class $A$ with $k$ marked points.  
Elements of this space are equivalence classes $[h, z_1,\dots, z_k]$ of elements $(h,
z_1,\dots, z_k)$, where $h:S^2\to P$ is a
somewhere injective (i.e. nonmultiply covered) 
$J$-holomorphic map and the $z_i$ are distinct points in $S^2$.  Elements
$(h, z_1,\dots, z_k)$ and $(h', z_1',\dots, z_k')$ are equivalent if there is $\ga\in
\PSL(2,\C)$ such that
$h' = h\circ \ga$, and $z_j = \ga(z_j')$ for all $j$.
                                                                                                                                                                                              
When $J$ is  generic, $\Mm_{0,k}(P,J,A)$ is a manifold of dimension 
$$
\dim(P) +
2c_1(TP)(A) +2k-6 = 2n + 4 + 2c_\phi(A),
$$
when $\dim P = 2n+2$, $k=2$ and $A$ is a section class of $P = P_\phi$. However it is
usually not compact.  Its compactification $\oMm_{0,k}(P,J,A)$ consists of 
$J$-holomorphic stable maps, i.e. of  equivalence classes $\tau = [\Si, h, z_1,\dots, z_k]$
of elements $(\Si, h, z_1,\dots, z_k)$.  Here the domain $\Si$  is a connected union of
components $\Si_i, i = 0,\dots, \ell-1,$ each of which has a given identification
with $S^2$.  (Note that we consider $\Si$ to be a topological space: the labelling
of its components is a convenience and not part of the data.)
The intersection pattern of the components can be described by a tree
graph with $\ell$ vertices, where each edge corresponds to an intersection
point of the components corresponding to its vertices.   No more
than two components meet at any point.  There are also $k$ marked points
$z_1,\dots, z_k$ placed anywhere on $\Si$ except at an intersection point of two
components.  The restriction $h_i$ of the map $h$ to $\Si_i$ is assumed to be
$J$-holomorphic and we impose the {\it stability condition}
that $h_i$ is nonconstant
(though possibly a multiple covering) unless $\Si_i$ contains at least  $3$ special
points.  (By definition, {\it special points} are either points of intersection  with other
components or marked points.)  Also $h_*([\Si]) = \sum_i (h_i)_*[\Si_i] = A$.   Finally,
we divide out by all holomorphic reparametrizations.  In other words if
$\ga:\Si\to\Si$ is a holomorphic map such that $ \ga(z_j') = z_j$ for all $j$ we consider
$(\Si, h\circ \ga, z_1',\dots, z_k')$ and $(\Si, h, z_1,\dots, z_k)$ to be
equivalent. 
Note that when $k=2$  the group of such reparametrizations for which
$z_j' = z_j$ for all $j$  has (real) dimension at
least $2$.  (The other reparametrizations  in some sense do not count since
they can be gotten rid of by normalizing the choice of $z_j$.)

This defines the elements $\tau$ of $\oMm=\oMm_{0,k}(P,J,A)$.  
We write $\ev$ for
the  evaluation map $$
\ev: \oMm_{0,k}(P,A,J)\to P^k:\quad [\Si, h,z_1,\dots,z_k]\to (h(z_1),\dots,h(z_k)).
$$

Each such $\tau$ has a finite automorphism group $\Ga_\tau = \Aut([\Si,h,
z_1,\dots, z_k])$.
Given a representative $(\Si, h, z_1,\dots, z_k)$  of $\tau$ we may identify $\Ga_\tau$
with the group  of all holomorphic maps $\ga:\Si\to\Si$
such that $h\circ \ga = h,$ and $ \ga(z_j) = z_j$ for all $j$.  
For most elements this group is trivial, and it is always finite because of the
stability condition.   It may be nontrivial if any of the $h_i$ are multiple
coverings, or if two different $h_i$ have the
same image. (Note that $\ga$ might permute the components of $\Si$.)
It is the presence of these automorphism groups that makes it impossible to 
find a manifold that regularizes $\oMm=\oMm_{0,k}(P,J,A)$.

The space of all stable maps has a natural topology in which $\oMm$ is compact.
We will later consider a small neighborhood $\Ww$ of $\oMm$ in this space. 
$\Ww$  has a natural {\it coarse} stratification in which each stratum 
consists of elements $\tau = [\Si_\tau, h_\tau, z_1,\dots, z_k]$ whose domain $\Si_\tau
= \cup_i\Si_i$ has fixed topological type, where the  marked points $z_1,\dots, z_k$ are
on specified components of $\Si$, and where the decomposition $A = \sum A_i$ is
given.  It also has a {\it fine} stratification  in which we add the requirement that  the
isomorphism class of $\Ga_\tau$ is fixed on each  stratum.  Each such stratum in
$\Ww$ has a smooth topology modelled on the space of smooth (i.e. in some Sobolev
class $L^{k,p}$) sections of a suitable bundle over the (fixed) domain $\Si$.  However, 
this smooth structure is not extended over $\Ww$ as a whole: see \S 4.
 It is not hard to check that the  index of the linearization $Dh_\tau$ of the
Cauchy-Riemann operator at $h_\tau$  depends only on the number $\ell$ of
components of the domain $\Si_\tau$ and equals  
\begin{equation}\label{eq:dim2}
{\rm ind}\, h_\tau = \dim \Mm - 2(\ell-1) = 2n + 2c_\phi(A) +  2k + 2 - 2\ell,
\end{equation}
where the second equality holds when $ P = P_\phi$ and $A$ is a section class.  

There is a special case in which $P$ is a single point.  In this case all maps are trivial,
and the resulting space $\oMm_{0,k}(pt,J,A)$ is usually written $\oMm_{0,k}$.  Its
elements are called {\em stable curves}.   The stability condition now says that every   
(nondegenerate) component of $\Si$ contains at least $3$ special points.  When $k < 3$
the space $\oMm_{0,k}$ reduces to a single point, that for consistency is taken to be
represented by the degenerate Riemann surface consisting of a single point.  The space
$\oMm_{0,k}$ also is a single point, but now this is represented by a sphere with $3$
marked points.   There is a forgetful map
$$
f:\oMm_{0,k}(P,J,A)\to \oMm_{0,k}: \quad [\Si, h,z_1,\dots,z_k]\mapsto
[\Si', z_1,\dots,z_k],
$$
where $[\Si', z_1,\dots,z_k]$ is the stabilization of the domain $[\Si, z_1,\dots,z_k]$ of
$h$ that is obtained by collapsing to a point each  component of $\Si$ that does not
contain at least $3$ special points.  (The domain $[\Si, z_1,\dots,z_k]$ is called a {\em
prestable curve}.)
Note that, for $k\ge 3$,
 $\oMm_{0,k}$ is a smooth oriented compact manifold of dimension
$2k-6$:  see [HS2] for example. 
\MS

We now discuss stable maps in the context of our fibered space $P\to S^2$.
As always, $\Om$ is a compatible symplectic form on $P$.

\begin{defn}\label{def:Jfib}\rm  We say that an 
almost complex structure $J$  on $P$ is {\em fibered} if  the
following conditions hold.   At each point $x\in P$ let $H_x$ denote the
$\Om$-orthogonal to the tangent space to the fiber.  Then $J$ must
 preserve the splitting $T_xP_\phi = T_x({\rm Fiber})\oplus H_x$,
be $\Om$-compatible and be such that the projection $p$ is
$(J,j)$-holomorphic for the standard complex  structure $j$ on $S^2$.
\end{defn}

The next lemma explains the structure of $J$-holomorphic stable maps
with two marked points when $J$ is
fibered.

\begin{lemma}\label{le:stembr1}  When $P = P_\phi$, $J$ is fibered and $\si$ is a
section class, each element $[\Si, h,z_0,z_1]$ of $\oMm_{0,2}(P,J,\si)$ contains one
component $(\Si_0, h_0)$, called the stem, that is a section of $p:P\to S^2$.  All other
components lie in the fibers of $P$ and may be divided into  a finite number of
connected pieces called branches, each branch lying in a different fiber. \end{lemma}
\proof{}  When $J$ is  fibered,
 any $J$-holomorphic curve $h:S^2\to P$ that
represents the class of a section  must intersect each fiber of $P$
exactly once transversally since $p\circ h: S^2\to S^2$ is a holomorphic map of degree
$1$.
Further, any $J$-holomorphic curve
$h:S^2\to P$ that represents a class $B$ such that $B\cdot [M] = 0$ must
lie entirely in a fiber, and no class with $B\cdot [M] < 0$ has a
$J$-holomorphic representative.  Therefore, given  a $J$-holomorphic
stable map $[\Si, h]$ in the class $\si$,  in the corresponding decomposition $\si =
\sum B_i$ of $\si$ (where $B_i = h_*(\Si_i)$)
 we must have $B_i\cdot [M] \ge 0$ for all $i$. 
 It follows
that $B_i\cdot [M] = 0$ for all $i$ except one (that we choose to be $i = 0$).  Each
other component  lies entirely in a fiber.  Moreover the union of the components
in a given fiber must be connected since $h(\Si)$ is connected 
and  the stem $h(\Si_0)$ intersects
each fiber just once.  The  result is now
immediate. \QED

  Each branch can be considered as an element of some moduli space
$\oMm_{0,k}(M,J_M,B)$, where $1\le k \le 3$.  For each branch contains  one 
marked point where it meets the stem and may or may not also contain
$z_1$ or $z_2$.

\subsubsection{Gromov--Witten invariants}\label{ssec:GW}

Siebert's point of view is the following.   In order to count the number of
isolated $J$-holomorphic $A$-curves of  genus $0$ in $(P,J)$ through the $k$ classes
$\al_1,\dots,\al_k \in H_*(P)$, one considers the space $\oMm_{0,k}(P,A,J)$ of
$J$-holomorphic   genus $0$ stable maps in class $A$ with $k$ marked points. 
This space  supports
 a {\it virtual fundamental class} 
$$
\Cc_{A,k}^{hol}(P,J) \in H_d^{BM}(\oMm_{0,k}(P,A,J), \Q)
$$
of dimension $d = \dim P + 2 c_1(A) + 2k - 6$, where $H^{BM}$ denotes the
Borel--Moore homology.   In the nicest cases, $\oMm_{0,k}(P,A,J)$ is a compact
manifold of dimension $d$ and  $\Cc_{A,k}^{hol}(P,J)$ is simply its fundamental
class.  More generally, (e.g. in the semi-positive case, see~[MS1])
$\oMm_{0,k}(P,A,J)$ is a compact stratified  space whose top stratum is a
manifold of dimension $d$ and all other strata  have dimension $\le  d-2$.  In
this case too  $\oMm_{0,k}(P,A,J)$ carries a fundamental class and 
$\Cc_{A,k}^{hol}(P,J) = \oMm_{0,k}(P,A,J)$.

Given $k$ homology classes
$v_1,\dots,v_k$ in $P$, the Gromov--Witten invariant $n_P(v_1,\dots, v_k; A)$ is
defined to be the intersection number 
$$
n_P(v_1,\dots, v_k; A) = \ev(\Cc_{A,k}^{hol}(P,J))\cdot_{P^k} (v_1\times\dots\times
v_k) 
$$
of the pushforward
$\ev_*(\Cc_{A,k}^{hol}(P,J))$ of $\Cc_{A,k}^{hol}(P,J)$ by the evaluation map
with the class $v_1\times\dots\times v_k$ in $H_*(P^k,\Q)$.  By definition, this
number is zero unless the dimensions of the cycles are complementary in $P^k$.
Thus if $P = P_\phi$ and $k=2$ we need
$$
d + \dim v_1 + \dim v_2  = 2(2n+2),
$$
which reduces to the dimension condition $2c_\phi(A) + \dim v_1 + \dim v_2 = 2n$
in equation~(\ref{eq:dim1}) above.  Note also that $n_P(v_1,\dots, v_k; A)
= 0$  if $A\ne 0$ and any $v_i = [P]$ since in that case it is
impossible for the intersections to be isolated.

The invariant $n_P(v_1,\dots, v_k; A)$ is the most simple-minded 
Gromov--Witten invariant.  More generally, one can consider the forgetful
map $$
f: \oMm_{0,k}(P,A,J)\to \oMm_{0,k},
$$
take a class $c\in H_*(\oMm_{0,k},\Q)$, and then set
$$
n_{P,c}(v_1,\dots, v_k; A) = 
(f\times \ev)_*(\Cc)\cdot (c\times v_1\times\dots\times
v_k), $$
where the intersection is taken in $\oMm_{0,k}\times P^k$.
Thus our previous definition corresponds to
taking $c$ to be the fundamental class of $\oMm_{0,k}$.

When $k\le 3$, $\oMm_{0,k}$ is a single point, and so there is only one invariant.
We also need to consider the case $k> 3$.  
Let  $\chi = [pt] \in H_0(\oMm_{0,k})$.
The invariant 
$$
 n_{P,\chi}(v_1,\dots, v_k; A) 
$$
intuitively speaking counts the number of $J$-holomorphic 
 maps  $h:S^2\to P$ such that
$h(z_i)\in v_i$ for a fixed set of points $z_1,\dots, z_k$.  In particular, if
$k=4$ fixing the points $z_i$ is equivalent to fixing their cross ratio.
Hence 
$n_{P,\chi}(v_1,\dots, v_4; A) $  is calculated using the class
$\Cc_{A,4}^{hol}(P,J)\cap [\oMm_t]$  in $H_*(\oMm_t)$, 
where $\oMm_t$ is the subset
of $\oMm = \oMm_{0,4}(P,J,A)$ on which the cross ratio 
of the marked points is fixed
at $t$.\footnote
{
The cross ratio on $\oMm_{0,4}(P,J,A)$ can be identified with the forgetful map
to the space $\oMm_{0,4} = S^2$ of genus zero stable curves with $4$ marked
points. For generalizations of this, see Hofer--Salamon~[HS2].}

 An important identity, that is equivalent to the associativity of quantum
multiplication, is the following:
\begin{eqnarray}\label{eq:assoc1}
n_{P,\chi}(v_1,\dots, v_4; A) = \sum_{A= A_1+A_2, \,\al} 
 n_{P}(v_1, v_2, e_\al'; A_1)\,\cdot\,
n_{P}(f_\al', v_3, v_4; A_2)
\end{eqnarray}
where $\{e_\al'\}_\al$ is a basis for the homology $H_*(P,\Q)$, and
$\{f_\be'\}_\be$ is the dual basis with respect to the intersection pairing. 
Further, if $v_4 = [P]$, we have
$$
n_{P,\chi}(v_1, v_2, v_3, [P]; A) = n_{P}(v_1, v_2, v_3; A).
$$
(This is consistent with~(\ref{eq:assoc1}) since in this case the only term that
 contributes to the sum  has $A_2 = 0$.)  Another important
point is that if $w\in H_{(\dim P) - 2}(P)$ then 
$$ 
n_P(v_1,v_2, w; A) = 
(w\cdot A)\, n_P(v_1,v_2; A). $$
 In particular, if $P= P_\phi$ and
$\si$ is a section class then
\begin{eqnarray}\label{eq:23}
n_{P}(v_1,v_2, [M], [M]; \si) = n_P(v_1,v_2,[M]; \si) = n_P(v_1,v_2; \si).
\end{eqnarray}
\SS

Finally, we observe that part (i) of Proposition~\ref{prop:GW}
 follows immediately from the definitions.  This states that
 $n_P(\io(a), \io(b), v; \io(A))$ is always $0$.  To see this, 
 observe that if $J$ is fibered Lemma~\ref{le:stembr1} implies that every
$J$-holomorphic stable map in class $\io(A)$ has image in a single fiber.  Hence if we
represent the classes $\io(a), \io(b)$ by cycles in different fibers the evaluation map
$$
ev:\oMm_{0,3}(P_\phi, J,\io(A)) \to P_\phi^3
$$
never intersects the representative of $\io(a)\times \io(b)\times v$.

It is also easy to check that the rest of Proposition~\ref{prop:GW} holds when $P$ 
is so nice that there is a fibered $J$ such that $\oMm_{0,k}(P,J,A)$ supports a
fundamental class whenever $A$ is a section or fiber class.   For then the relevant Gromov--Witten invariants can
be calculated by simply counting $J$-holomorphic curves.\footnote{
This is the case considered
in [Sd1] and  [LMP2].  To complete the argument in this case one must
use  some elementary  transversality arguments that are similar
to those in Lemma~\ref{le:transv} below.}

\subsection{Proofs of Propositions~2.2 and~2.5.}

In this section, we will prove these propositions modulo  two statements
Claim~\ref{cl} and Lemma~\ref{le:xx} about the behavior of
Gromov--Witten invariants.  To simplify our notation we will usually omit
mention of the inclusion map $\io: H_*(M)\to H_*(P)$, instead using the
convention that $a,b,c, A,B$ denote elements of $H_*(M)$ and $v,w,\si$
denote elements of $H_*(P)$.

 \subsubsection{Calculating $\Psi_{\phi,\si}$ for trivial
bundles}

 Proposition~\ref{prop:const}  says that $\Psi_{\phi,\si_0}$ is the
identity when $\phi$ is the constant loop $*$ and $\si_0$ is the flat section.
Thus we must show:

\begin{lemma} {\rm (i)} $n_P(a,b;\si_0+ B) = 0$ unless $B=0$, \SS

\NI
{\rm (ii)}
$n_P(a,b;\si_0) = a\cdot_M b$.  
\end{lemma}

\NI
{\bf Proof of (ii)}.

Let $\Om$ be a product form on $P_* = M\times S^2$, $\al,\be$ be cycles in
different fibers of $P_*$ and $J$ be a product   almost complex structure
$J_M\times j$  on $P_*$.   
Then the projections
$P_*\to M,$ $ P_*\to S^2$ are holomorphic, so that any $J$-holomorphic
curve $\tilde f: S^2\to M\times S^2$ in class $\si_0$ projects to a constant in
$M$.   Thus,  the space $\Mm_{0,2}(P,\si_0,J)$ of
unparametrized $J$-holomorphic spheres in class $\si_0$ is a compact manifold
 that can be identified with $M$.  In
this case the virtual fundamental class  $
\Cc_{A,k}^{hol}(P,J)$ is simply the fundamental class of $\Mm_{0,2}(P_*,\si_0,J)$.
It follows immediately that $n_P(a,b;\si_0) = a\cdot_M b$.\QED

\NI
{\bf Sketch Proof of (i)}.

Now consider calculating $n_P(a,b;\si_0+ B)$ where $B\ne 0$.
First observe that when $J = J_M\times j$ there is a well defined map
$$
gr: \oMm_{0,3}(M,B,J_M) \to \oMm_{0,2}(M\times S^2, \si_0 + B, J)
$$
given on $\tau = [\Si_\tau, h, z_0, z_1, z_2]$ as follows.  If the domain
$\Si_\tau$ has just one component, then 
there is a unique map $\ga: \Si_\tau\to S^2$ that takes the $3$ marked
points to $q_0 = 0, q_1 = 1, q_2= \infty$.  Then set
$gr(\tau) = [\Si_\tau, \widehat{h}, z_0, z_1]$ where
$$
\widehat{h} (z) = (h(z), \ga(z))\in M\times S^2.
$$
Thus, if we identify $\Si$ with $S^2$ via $\ga$, $\widehat h$ is the graph
of $h$.

In general, let us define the
projection $p_j(z_i)$ of the marked point $z_i$ onto the component $\Si_j$ as
follows. If $z_i\in \Si_j$ then $p_j(z_i) = z_i$.   Otherwise it is the unique point
on $\Si_j$ where the connected component of $\Si_\tau - \Si_j$ containing
$z_i$ is attached to $\Si_j$.   It is
easy to see that for each domain $[\Si_\tau,  z_0, z_1, z_2]$ there is
precisely one component  onto which the projections of the three marked
points are distinct.  Let us call this component $\Si_0$.  It contains 
three distinct points $ p_0(z_i)$ and so we can define $\widehat{h}$ on this
component as before.  This is the stem of $gr(\tau)$.  The other
components are lifted as branch components in the obvious way.

It is not hard to check that the map $gr$ is continuous.  Moreover, its
image  is precisely the set of  elements of $\oMm_{0,2}(M\times S^2,
\si_0 + B, J)$ that take $z_i$ into $ M\times q_i$ for $i= 0,1$.  Hence,
$gr$ is essentially surjective. (More precisely, its image contains
 all stable maps, but the marked points might be in the wrong place.)
  On the other hand, when we calculate the
Gromov--Witten invariant $n_P(a,b;\si_0 + B)$  it is all right to
have the marked points in designated fibers.
Indeed,
let us consider the evaluation map 
$$
\ev_{01}: {\rm Im\,} gr\to M\times M,
$$
where  the range  $M\times M$ 
  is identified with
$(M\times \{q_0\})\times (M\times \{q_1\})$ in $P\times P$.
Observe that the formal codimension of its image is the same as that of the map
$$
ev: \oMm_{0,2}(M\times S^2,\si_0 + B, J) \to P\times P.
$$
We will prove an equivalent of the following claim in Lemma~\ref{cl1}.

\begin{claim}\label{cl}  If $\Cc$ is the virtual
fundamental class  of $ \oMm_{0,3}(M,B,J_M)$,  then the invariant
$n_P(a,b;\si_0 + B)$ is the intersection number of 
$$
(\ev_{01}\circ gr)_* (\Cc)
$$
 with the class $a\times b$ in $M\times M$
\end{claim}

Now observe that the map $ev_{01}\circ gr$ reduces dimension:
given 
$$
\tau = [\Si, h, z_0, z_1, z_2]\in \oMm_{0,3}(M,B,J_M),
$$
the element $gr(\tau)$ varies as $z_2$ moves, but $ \ev_{01} (gr(\tau))$
remains unchanged.  It follows easily that, for dimensional reasons, $a\times b$
can be perturbed in $M\times M$ to be disjoint from a representing cycle for
the class  $(\ev_{01}\circ gr)_* (\Cc)$.    Hence $n_P(a,b;\si_0 + B)=0$.

\subsubsection{Proof of the composition rule}\label{ssec:comp}

It follows immediately from the definition of $\Psi$ 
(see  Definition~\ref{def:psi}) and the
remarks  at the end of \S~\ref{ssec:GW} that, if ${e_i}$ is a basis for
$H_*(M)$ with dual basis $\{f_j\}$, then 
$$
\Psi_{\phi,\si}(a) = \sum_{B,i} n_{P_\phi}(a,[M],e_i);\si + B) f_i\otimes
e^{-B}. 
$$
We would like to use equations~(\ref{eq:assoc1}) and~(\ref{eq:23})
in \S2.2.2  to
make the following calculation: 
\begin{eqnarray*}
\Psi_{\psi,\si'} (\Psi_{\phi,\si}(a))  & = & \sum_{B,i}
 n_{P_\phi}(a,[M], e_i;\si + B) \Psi_{\psi,\si'} (f_i)\otimes e^{-B}\\
& = & \sum_{B, B'}\sum_{i,k}
 n_{P_\phi}(a,[M], e_i;\si + B)\,\cdot\\
 & & \qquad\qquad\quad n_{P_\psi}(f_i,[M],e_k;\si' + B')
f_k\otimes e^{-B-B'}\\ & = & \sum_{A,k} n_{P_{\psi*\phi},\chi}(a,[M],
[M],e_k, \si\#\si' + A)\otimes e^{-A}\\\vspace{.1in}
& = & \Psi_{\psi*\phi,\si\#\si' }(a).
\end{eqnarray*}

The only problem here is that we are considering Gromov--Witten invariants
in different spaces.  To overcome this difficulty we must define a space $\Xx$ that
contains all the spaces $P_\phi, P_\psi$ and $P_{\psi*\phi}.$
To this end,  let $\pi: (\Ss,j_\Ss) \to
\De$ be a singular holomorphic fibration over the unit disc in $\C$, with generic fiber
$\Si_t= S^2$ and singular fiber $\Si_0$ over $\{0\}$ formed by the union of two
transversally intersecting $2$-spheres.  Thus there are local complex coordinates
$(z,w)$ round the singular point  in $\Si_0$ such that
$$
\pi(z,w) = zw.
$$
Observe that each component of $\Si_0$ has self-intersection
$-1$.  Choose 
four disjoint smooth sections $s_i:\De\to \Ss$ such that $s_1(0), s_2(0)$ lie in one
component of $\Si_0$ and $s_3(0), s_4(0)$ lie in  the other, and so
that the cross ratio of the four points $s_i(t)$ on $\Si_t$ uniquely determines
$t$ when $t\ne 0$.  This is a holomorphic model for the degeneration of the
complex structure of the sphere $\Si_t$ with $4$ marked points as a bubble is
formed that contains two of the marked points.  

We next construct a smooth fibration $q:\Xx\to \Ss$ with fiber $M$
such that
$$
X_t = \Xx|_{\Si_t} \cong P_{\psi}\# P_{\phi},\quad X_0=\Xx|_{\Si_0} = P_\phi\cup
P_\psi, 
$$
where in the union $P_\phi\cup P_\psi$ the fibers
$M_{\psi,0}$ are
$M_{\phi,\infty}$ are identified.
This can be done very explicitly using the description of 
$P_{\psi}\# P_{\phi}$ in \S 1.  Moreover, we can suppose that $\Xx$ has a
symplectic form $\Tilde\Om$ and a compatible $\Tilde J$ that restrict
appropriately on each $X_t$.  (In particular, the projection $\Xx\to \Ss$ is
holomorphic, which implies that  $\Tilde J$ restricts to a fibered almost complex
structure $J_t$ on each $X_t$.) 

Now   let $\Tilde\al_i, i = 1,2,3,4$ be relative cycles in $(\Xx, \p\Xx)$ whose images by
$q:\Xx\to \Ss$ equal the images of the sections $s_i:\De\to \Ss$ and are such that:
$$
\Tilde \al_1\cap X_0 = \al, \quad \Tilde \al_i = q^{-1}({\rm Im\,}s_i),\;i =
2,3,\quad \Tilde \al_4\cap X_0 = \be,
$$
where $\al,\be$  represent the classes $a,b$. 
Note that every $\Tilde J$-holomorphic stable map $\tau$  that 
meets these cycles $\Tilde \al_i$
must lie entirely in one of the sets $X_t$ (since the projection $\Xx\to \De$ is
holomorphic). Moreover, the value of $t$ is determined by the cross ratio of the
marked points on $\tau$.

Let $\Tilde a_i\in
H_*(\Xx,\p\Xx)$ be the classes represented by the $\Tilde\al_i$.
Then, even though $\Xx$ has a boundary, one should be able to make sense of  
Gromov--Witten invariants on $\Xx$.  It follows easily from the  calculation given
above that Proposition~\ref{prop:comp} would hold if we could prove the
following lemma.

\begin{lemma}\label{le:xx}  {\rm (i)}  If $\si$ is a section class in $X_t =
P_{\psi*\phi}$ and $\Tilde \si$ is its image in $H_2(\Xx)$, for each  $t\ne 0$ we have: 
$$
n_{\Xx,\chi}(\Tilde a_1,\dots, \Tilde a_4; \Tilde \si) = n_{X_t}(a, [M], [M], b; \si).
$$

\NI
{\rm (ii)} The above statement also holds  when $t=0$ if the right hand
side is defined by:
$$
n_{X_0}(a,[M], [M], b; \si) = \sum_{\si = \si_1\#\si_2,\, i
}n_{P_\phi}(a,[M],e_i;\si_1)\; n_{P_\psi}(f_i,[M],b;\si_2). $$
\end{lemma}

To justify (i) note that $n_{\Xx,\chi}$ is calculated by restricting attention 
to the part
of the space $\oMm_{0,4}(\Xx,\Tilde J, \Tilde \si)$ on which the cross ratio is
fixed, and, as pointed out above, the only elements in $\oMm_{0,4}(\Xx,\Tilde
J, \Tilde \si)$ that meet the cycles $\Tilde\al_i$ lie in  a
single $X_t$.   

A proof of this lemma will be given in \S4.3.3 below.

\section{Proofs of the other properties of $P_\phi$}

In this section we complete the proofs of the results 
announced in \S1.1 and \S 1.2,
assuming that Theorem~\ref{main} is proven and 
that Propositions~\ref{prop:GW}, \ref{prop:const} and~\ref{prop:comp}
hold.  Proposition~\ref{prop:GW} itself is proved in \S4.3.2.

\subsection{The homomorphism $\rho$.}

In order to use the maps
$\Psi_{\phi,\si}$ to define $\rho$
we must  make a canonical choice of section $\si_\phi$ that (up to  equivalence)
satisfies the composition rule 
$$
\si_{\psi*\phi} = \si_\phi\#\si_\psi,
$$
where $\si_\phi\#\si_\psi$ denotes the  obvious union of the sections
$\si_\phi, \si_\psi$ in
$P_{\psi*\phi} = P_\psi\#P_\phi$.  Unfortunately, it is not always possible to
do this if one just considers ordinary sections.  To get around this difficulty, we
define the notion of a generalized section or $\R$-section.

\begin{defn}\rm We will say that $\si$ is an $\R$-section of
$P_\phi$ if it is a finite sum $\sum\la_i\si_i, \la_i\in \R,$ of sections such that
$p_*([\si])$ is the fundamental class of $S^2$.  Two $\R$-sections $\si,\si'$ are
equivalent if  
$$
u_\phi(\si)= u_\phi(\si'),\quad c_\phi(\si)= c_\phi(\si').
$$
 \end{defn}

The following lemma is immediate. (Part (iii) follows from Lemma~\ref{le:glsect}.)

\begin{lemma}\label{le:si} {\rm (i)} If the classes $[\om]$ and $c_1(TM)$ are linearly
independent on $H_2^S(M)$, there is for each $\phi$  a unique 
equivalence class of $\R$-sections 
$\si_\phi$ such that  $$
 u_\phi(\si_\phi) =
0,\quad c_\phi(\si_\phi) = 0. 
$$ 
\NI {\rm (ii)}  If the classes $[\om]$ and $c_1(TM)$ are linearly
dependent on $H_2^S(M)$ and $[\om]\ne 0$, we choose 
$\si_\phi$ to be the
unique  equivalence class of $\R$-sections 
 such that  $
 u_\phi(\si_\phi) = 0
$.  If $[\om] = 0$ on $H_2^S(M)$ but $c_\phi\ne 0$, we choose 
$\si_\phi$ so that  $
c_\phi(\si_\phi) = 0$.  If both $[\om]$ and $c$ vanish on $H_2^S(M)$,
then there is only one section class up to  equivalence, and we take that for $\si_\phi$.

\NI
{\rm (iii)}
In all cases, $\si_{\psi*\phi} = \si_\phi\#\si_\psi.$
\end{lemma}

\begin{rmk}\rm 
Note that it may happen that the classes $[\om]$ and $c_1(TM)$ are linearly
dependent on $H_2^S(M)$ while their extensions $u_\phi$ and $c_\phi$ are not
linearly dependent on $H_2(P_\phi)$: see Example~\ref{ex:PD}.   Therefore, there
may be no $\R$-section such that
 $ u_\phi(\si_\phi) =
0,$ $ c_\phi(\si_\phi) = 0. $   This is why we used  weaker conditions
in (ii) above.  Note  however that  when $\om
= 0$ on $H_2^S(M)$  it follows from the proof of Proposition~\ref{prop:ring2}
in \S\ref{ssec:nonsq} that the section $\si_\phi$ defined as above is
integral and that $u_\phi(\si_\phi) = c_\phi(\si_\phi) = 0$.
 \end{rmk} 

One should think of the $\R$-section $\si_\phi$ as an average of the sections in
$P_\phi$.  The effect of enlarging the Novikov ring to $\La_R$ is 
thus to make enough room
to take this average.  More precisely,  observe that the definition of the map
$\Psi_{\phi,\si_\phi}$ still makes perfect sense provided that one allows the
coefficients $B$ to belong to $\Hh_R = H_2^S(M,\R)/\!\!\sim$ so that the sum
$\si_\phi + \io(B)$ can be integral.  In the next proposition 
$\Aut_{\La_R}(QH_*(M,\La_R))$ denotes the group of  
automorphisms of the $\La_R$-module
$QH_*(M,\La_R)$.

\begin{prop}\label{trho}  The map
$$
\rho':\; \pi_1(\Ham(M,\om))\to \Aut_{\La_R}(QH_*(M,\La_R))
$$
defined by $\rho'(\phi) = \Psi_{\phi,\si_\phi}$ is a group homomorphism.
\end{prop}
\proof{}  This follows immediately from 
Propositions~\ref{prop:const} and~\ref{prop:comp} and Lemma~\ref{le:si}.\QED

Following Seidel, we now want to show that
the automorphisms $\rho_\phi' = \rho'(\phi)$
 respect the structure of $QH_*(M)$ as
a right $QH_*(M)$-module, i.e.
$$
\rho_\phi' (a*_M b) = (\rho_\phi' a)*_M b.
$$
To this end, it is
convenient to introduce the following more flexible notation
for the quantum product.
Let $\{e_j\}$ be a basis for $H_*(M)$ with dual basis 
$\{f_j\}$.  Then, for any section $s':H_*(M)\to H_{*+2}(P)$ 
$$
\io(e_i)\cdot_P s'(f_j) = e_i\cdot_M f_j = s'(e_i)\cdot_P \io(f_j) =
\de_{ij}. 
$$
The following lemma shows that we can modify $s'$ to a section $s$ such that
in addition $s(e_i)\cdot s(e_j) = s(f_i)\cdot s(f_j) =0$ for all $i,j$.
To simplify our notation we will from now on omit mention of the inclusion
$\io: H_*(M)\to H_*(P)$, but will keep the convention that the letters
$a,b,c,e,f$ refer to homology classes on $M$ and $u,v,w$ to those in $P$.

\begin{lemma}\label{le:basis}  In the above situation, let $d(e)$ denote the
dimension of the class $e\in H_*(M)$, let $2n = \dim(M)$, set $q_{ij} =
s'(e_i)\cdot s'(e_j)$, and define $$
s(e_i) = \left\{\begin{array}{lll} s'(e_i) & \mbox{if} & d(e_i) < n-1\\
s'(e_i) - \frac 12\sum_j (-1)^{n-1}q_{ij}\, f_j & \mbox{if} &
d(e_i) = n-1\\ s'(e_i) -  \sum_j (-1)^{d(e_j)}q_{ij} \,f_j &\mbox{if} &
d(e_i) > n-1, \end{array}\right.
$$
where the sums are over those $j$ for which $f_j$ has the appropriate
dimension, i.e. $d(f_j) + d(e_j) = 2n$.  Then, if the $s(f_j)$ are defined by
linearity, $$
e_i\cdot_P s(f_j)  =\de_{ij},\quad s(e_i)\cdot_P s(f_j) = 0
$$
for all $i,j$.
\end{lemma}

\proof{} This is an easy calculation that is left to the reader.\QED

From now on we will assume that the $e_j, f_j$ and $s$ are as above.
Then the elements $\{e_j, s(e_k)\}$ form a basis for $H_*(P)$ whose dual
basis is given by 
$\{s(f_j), $ $ f_k\}$. Similarly, the basis $\{f_j, s(e_k)\}$ has dual basis
$$
\{(-1)^{\eps_j} s(e_j), f_k\},\quad\mbox{where}\;\;\eps_j \equiv d(e_j) \pmod
2.
$$

Observe that in our present notation
\begin{eqnarray*}
a*_M b & = &  \sum_{j, B\in \Hh_R}\;
n_M(a,b,e_j;B) f_j \otimes e^{-B} \in QH_*(M,\La), \\
\Psi_{\phi,\si} (a) & = & \sum_{j,B\in \Hh_R} n_P(a, [M],e_j; \si + B)
\, f_j\otimes e^{-B}.
\end{eqnarray*}
It follows easily that the triple invariants $n_P(a, b, c\,;\si)$
 are the coefficients of $\Psi_{\phi,\si}
(a)*_M b$. Thus, we have:

\begin{cor}\label{cor:GW1}
$$
\Psi_{\phi,\si} (a)*_M b = \sum_{B, j} n_P(a,b, e_j; \si + B)\, f_j\otimes e^{-B}.
$$
\end{cor}

 Similarly, the horizontal quantum product $*_{H,\si}$
that was introduced in \S 1.3 may be written as:
\begin{eqnarray*}
u*_{H,\si} v & = &  \sum_{j,B\in \Hh_R}
 n_P(u,v, e_j; \si+B) \,s(f_j) \otimes e^{-B}
 \\ & & \qquad\qquad +\; n_P(u,v, s(e_j); \si+B)\, f_j \otimes e^{-B}. 
\end{eqnarray*}

\begin{prop}\label{prop:qsi} 
 {\rm (i)}  $\Psi_{\phi,\si}(a) = ([M]*_{H,\si} a)\cap [M]$.\SS

\NI
{\rm (ii)} Let $Q_\si = \Psi_{\phi,\si}([M]) =
([M]*_{H,\si}[M])\cap [M]$. Then
$$
\Psi_{\phi,\si}(a) = Q_\si *_M a.
$$
\end{prop}
\proof{} The first statement follows immediately from 
Definition~\ref{def:psi}
since $s(f_j)\cap[M] = f_j,\,  f_j\cap [M] = 0$.  To prove (ii)
observe first that
\begin{eqnarray*}
\Psi_{\phi,\si}(a) & = & \left([M]*_{H,\si} a\right)\cap [M]\\
& = & \sum_{j,B\in \Hh_R}
 n_P([M], a, e_j; \si+B)\; f_j \otimes
e^{-B}\\
& =&\sum_{B,A,k,j} n_P([M], [M], e_k; \si +B-A)
\,\cdot\, n_M(f_k, a,e_j;A)\otimes e^{-B}
\end{eqnarray*} 
by Proposition~\ref{prop:GW}(iv).
But
\begin{eqnarray*}
Q_\si *_M a & = &  \sum_{C,k} n_P([M], [M], e_k; \si +C)\; f_k*_M a
\otimes e^{-C}\\
 & = & \sum_{C, A, k,j} n_P([M], [M], e_k; \si + C)\,\cdot\, 
n_M(f_k,a,e_j;A) \;f_j \otimes e^{-A-C}.
\end{eqnarray*}
The result follows.\QED

\begin{cor}\label{cor:qsi} 
$\Psi_{\phi,\si} (a*_M b) = (\Psi_{\phi,\si}(a))*_M b.$
\end{cor}

It follows from the above proposition that, for all sections 
$\si$,  $Q_\si =
\Psi_{\phi,\si}([M])$
 is a unit in
the ring $QH_{ev}(M,\La_R)$.
We now define the map  
$$
\rho:\pi_1(\Ham(M,\om))\to QH_*(M,\La_R)^\times
$$
by: 
$$
 \rho(\phi) = Q_{\si_\phi}.
$$

 \NI
{\bf Proof of Theorem~\ref{repr}}

We have to show that $\rho$ is a homomorphism.  
Proposition~\ref{prop:qsi} implies that
$$
\rho'(\phi) (a) = \rho(\phi)*_M a,
$$
and by  Proposition~\ref{trho} $\rho'$ is a
homomorphism. Hence result.\QED

The next result will be useful in \S3.3.  

\begin{lemma}\label{le:id}$\rho(\phi) = \mu \1\otimes e^{-A}$ if and only if
$$\begin{array}{lcll} 
n_P(a,b,[M];\si_\phi + B) &= & 0, & B\ne A\\
&= & \mu\,a\cdot b,& B = A.\end{array}
$$
\end{lemma}
\proof{}  This is immediate from the definitions.\QED

Following Seidel let us write
$$
QH^+ = \bigoplus_{k<n} H_{2k}\otimes\La_R,\quad\mbox{where }\; 2n = \dim M.
$$
The next result follows immediately from the definitions.

\begin{lemma}\label{le:unit}
{\rm (i)}  If the Gromov--Witten invariants
$n_M(a,b,c\,; B)$ vanish whenever $c_1(B) > 0$ then $QH^+ *_M QH^+\subset QH^+$.
\SS

\NI
{\rm (ii)} If $QH^+ *_M QH^+\subset QH^+$ then every unit in
$QH_{ev}^\times$ has the form $\1 \otimes \la + x$ for some unit $\la \in \La_R$
and some $x\in QH^+$. \end{lemma}

Units in $QH_*(P, \La_R)$ of the form $\1 \otimes \la + x$ are in some sense
trivial. Observe that if a unit $\la = \sum_B\la_B e^{B}\in \La_R$ 
has more than one nonzero coefficient $\la_B$ then either $\la$ itself or its
inverse $\la^{-1}$ has infinitely many nonzero coefficients. It is not clear
whether such units can occur as values of $\rho$.  If so, it would mean that
$P_\phi$ had  nonvanishing  Gromov--Witten invariants for infinitely many
section classes.  

\begin{lemma}  If $\rho(\phi)= \1 \otimes \la + x$ where $\la\ne 0$ then
$\oI_c(\phi) = 0$. \end{lemma}
\proof{}  The hypothesis means that one of the numbers
$n_P([M], [M], pt;\si)$ is nonzero, which is possible only if $c_\phi(\si) = 0$.
Thus there is an (integral) section of $P\to S^2$ on which $c_\phi = 0$, which
immediately gives the result.\QED

In [Sd1] \S11 Seidel calculated $\Psi_{\phi,\si}$ in some cases, for example,  when
$\phi$ is rotation in the first coordinate of ${\CP}^n$, and when $\phi$ is 
rotation in the fiber of the rational ruled surface ${\bf F}_2$. 

The next  example is slightly easier because $H_2^S(M)$ has rank $1$.  We will
also calculate the homomorphisms $\oI_c, \oI_u$.

\begin{example}\label{ex:PD} \rm Let $M$ be total space of the nontrivial
$S^2$-bundle over $T^2$  and
let $\phi$ be the loop that rotates once around in the fibers
fixing the points
of the sections $T_{\pm}$ of self-intersection $\pm 1$.  Then
 $P_\phi$  fibers over $Z = T^2\times S^2$ with fiber $F = S^2$, and $M\subset
P_\phi$  can be identified with
the inverse image of $T^2\times pt$.  
In fact $P$ can be can be thought of as the
projectivization of $L\oplus \C$ where $L$ is a line bundle over $Z$ with
$ c_1(L) = \la + \mu,$ where 
$\la,\mu$ are the obvious 
generators of $H^2(Z)$ with $\la(T^2) = 1, \mu(S^2) = 1$.
We will use the same letters $\la,\mu$ for
the pullbacks of these classes to $P$.  There is another good class in $H^2(P)$
namely the vertical Chern class $\nu$ of the fibration $P\to Z$. 
 This is $2$ on
the fiber $F$  and takes the values $\pm 1$ on $S_\pm, T_\pm$, where
 $T_\pm, S_\pm$ are  copies  
 of  $T^2\times pt$ and $pt\times S^2$ in the two obvious 
sections  $Z_\pm$ of $P\to Z$.

The first claim is that the even part $H^{ev}(P)$ of $H^*(P)$ is generated by
$\la,\mu,\nu$ with the relations $\la^2 = \mu^2 = 0, \nu^2 = 2\la\mu $.  To see
this, note that $\PD(\nu) = Z_+ + Z_-$.  Hence $$
\PD(\nu^2) = (Z_+ + Z_-)^2 = (S_+ + T_+) - (S_- + T_-) = 2F.
$$
Also $\PD(\mu) = M$ and $\PD(\la) = Z_+ - Z_- - M$.   To check:
\begin{eqnarray*}
(Z_+ - Z_- - M)^2 &= &Z_+^2 + Z_-^2 -2M(Z_+ - Z_-) \\
&  = &  (T_+ + S_+) - (T_- + S_-) -
2(T_+ - T_-)\\
& = &  2F - 2F = 0.
\end{eqnarray*}
Also $\PD(\la\mu) = M(Z_+ - Z_- - M) = T_+ - T_- = F$.  (This is another check,
since obviously $\PD(\la\mu) = F$.)

Let us calculate $\oI_c(\phi) = c_\phi(\si)\pmod N$.  In our case the minimal 
(spherical) Chern number of $M$ $=2$.  Further the classes $S_\pm$ are
section classes of $P\to S^2$.  Hence 
$$
\oI_c(\phi) = 1\in \Z/2\Z.
$$
Next, identify $H_2(M,\R)/H_2^S(M,\R)$ with $ \R T$ where $T=T_+= T_-$.
  Observe  that 
$$
\PD(\la\nu) = (Z_+ - Z_- - M)(Z_+ + Z_-) = (T_+ + S_+ + T_- + S_- ) - T_+ - T_- = S_+ +
S_-
$$
 is spherical.  So are $\PD(\nu^2)$ and $\PD(\la\mu)$, while $\PD(\mu\nu) = T_+ + T_-$ is
not. Let us suppose that the symplectic form $\om$ on $M$ is such that
$$
\om(T_-) = \ka,\quad \om(F) = 2.
$$
Then  $u_\phi = (1+\ka) \la +  \eps\mu  +  \nu$ where $\eps$ is chosen so that
$$
0 = u_\phi^3 =  c^3(3(1+\ka)\eps + 1).
$$
Hence 
$$
\oI_u(\phi) = 2\eps\PD(\mu\nu) = 4\eps T= \frac {-4}{3(1+\ka)} T.
$$

Finally, let us calculate $\rho(\phi)$.  The first step is to calculate
$QH_{ev}(M)$.  Since the only
spherical class in $H_2(M)$ is $F$, it is easy to check that the quantum
products of the basis elements $[M] = \1, F, T_-, [pt]$ equal their usual
intersection product except in the following cases: 
$$ 
T_-* T_- = -[pt] + \1\otimes e^{-F},\quad [pt] * T_- = F\otimes e^{-F}.
$$
Next, we need to find all
Gromov--Witten invariants 
$$
n_P([M], [M], a; \si) = n_P([M], a; \si)
$$
 where
$\si$ is a section class.  Observe  there is a complex structure $J$ on $P$ for
which the submanifolds $Z_\pm$, $T_\pm, S_\pm$ are holomorphic, as are
the projections $p_Z: P\to Z$, $p: P\to S^2$.
 Therefore, any $J$-holomorphic section $C$
of $P\to S^2$ projects to a holomorphic sphere $\overline C$  in $Z$ and hence
lies in the ruled surface $p_Z^{-1}(\overline C)$.  It follows that  $C$ must lie in
a section class $\si = S_- + k F$ with $k\ge 0$.  In particular, it is regular. 
Hence the moduli space $\oMm_{0,2}(P, J, \si)$ carries a fundamental class of
the right dimension, and one can calculate $n_P([M], a; \si)$ by naively
counting curves.  In particular, the set of $J$-holomorphic curves in class
$S_-$ fills out $Z_-$ and so
$$ n_P([M], F; S_-) = 1,\quad n_P([M], T_-; S_-)
= -1,\quad n_P([M], a; S_+) = 0, $$
where the last equality holds for dimensional reasons: $n_P([M], a; \si) \ne
0$ only if $c_\phi(\si) \le 0$ by equation~(\ref{eq:dim1}) in \S2.1.
Hence
$$
\Psi_{\phi, S_-}[M] = T_-.
$$
(This also follows by applying [Sd1] Proposition 7.11.) 

Now observe that $M$ is spherically monotone,
and $\Hh = \Z$ is generated by $F$.  The class $\si_\phi$ must have the form $S_- +
\de F$ where, by Lemma~\ref{le:si}(ii) $u_\phi(\si_\phi) = 0$.  Hence
$$
\si_\phi = S_- + \de F,\quad \mbox{where} \quad \de = \frac{4+3\ka}{6+6\ka}.
$$
Thus
$$
\rho(\phi) = T_-\otimes e^{\de F}.
$$
 This is a unit in
$QH_{ev}(M, \La_R)$ of infinite order.  Hence $\phi$ itself must have infinite
order.
Observe finally that $c_\phi(\si_\phi)\ne 0$.  This is no contradiction since we are in
a degenerate case of Lemma~\ref{le:si}.\QED
\end{example}

\subsection{The vertical quantum cohomology of $P$}

Additively 
$QH_V^*(P, \La_R^*) = H^*(P;\La_R^*)$, where $\La_R^*$ is the Novikov ring of
$\Hh_R = \Hh_R(M)$ that satisfies the  finiteness condition dual to that
defining $\La_R$. Thus the elements of $\La_R^*$ have the form $\sum_{B\in
\Hh_R} \la_B^* e^B$ where, for each $\ka > 0$,
there are only finitely many nonzero coefficients $\la_B^*$ with $\om(B) <
\ka$.  The degree of $e^B$ is defined to be $2c_1(B)$ as before.  

Recall that the quantum cup product on $M$ is  the Poincar\'e dual to the
quantum intersection product on $M$.   In detail:  there
 is a dual pairing 
$$
QH^*(M,\La_R^*)\otimes QH_*(M, \La_R) \to  \Q:\;
\langle \sum\al_A^*
\otimes e^A,\sum a_B \otimes e^{-B}\rangle  =  \sum_{B} \langle
\al_B^*, a_{B}\rangle,
$$
and Poincar\'e duality $\PD = \PD_M$ is given by:
$$
\PD: QH^*(M, \La_R^*)\to QH_*(M, \La_R):\quad \sum \al_A\otimes e^A\mapsto
\sum \PD(\al_A) \otimes e^{-A}.
$$
The quantum cup product on $M$ is defined by
$$
\al*\be = \sum_B (\al *\be)_B \otimes e^B,
$$
where
$$
\langle (\al* \be)_B,\, c\rangle = n_M(\PD(\al), \PD(\be), c\,; B),
$$
and it is not hard to check that
$$
\PD(\al *\be) = \PD(\al) * \PD(\be).
$$
The vertical quantum product $*_V$ on
$QH^*(P,\La_R^*)$ is defined similarly.

 Theorem~\ref{main} implies that the
 Wang sequence in rational cohomology decomposes as a family of exact
sequences:
\begin{equation}\label{eq:wang2}
0\to H^{*-2}(M)\stackrel {\io_!}{\to} H^*(P_\phi) 
\stackrel{\io^*}\to
H^{*}(M)\to 0,
\end{equation}
where $\io_!: H^{*-2}(M)\to H^*(P)$ is the cohomology transfer:
$$
\io_!(\al) = \PD_P(\io_*(\PD_M(\al))).
$$
Observe that $\io^*$ is a ring homomorphism while $\io_!$ satisfies the
identity: $$
\io_!(\al\cup\be) = s(\al)\cup \io_!(\be),
$$
where $s: H^*(M)\to H^*(P)$ is any right inverse to $\io^*$, i.e. any map such
that $\io^*\circ s = {\rm id}$.  Here is the quantum analog.

\begin{prop}\label{prop:cohW} There is an exact sequence
$$
0\to QH^{*-2}(M, \La_R^*)\stackrel {\io_!}{\to} QH_V^*(P_\phi, \La_R^*) 
\stackrel{\io^*}\to
QH^{*}(M, \La_R^*)\to 0,
$$
of $\La_R^*$-modules where $\io^*$ is a ring homomorphism with respect to
the quantum products and $\io_!$ satisfies the identity
$$
\io_!(\al*_M\be) = s(\al)*_V \io_!(\be),
$$
for any additive splitting $s$ as above.
\end{prop}
\proof{}  This follows  by dualizing Proposition~\ref{prop:wang}.\QED

\begin{rmk}\label{rmk:sei} \rm
This proposition continues to hold whenever $P\to S^d$ is a fibration with
structural group $\Ham(M)$.  In fact, it is shown in [LMP3] that the
cohomology of $P$ is still additively split in this case. It is not hard to
check that all that is then  needed to establish the existence of the vertical
quantum multiplication is the  fiberwise symplectic structure of $P$.  Seidel
notes in [Sd2]  that the Wang
sequence~(\ref{eq:wang2}) with base $S^d$ can be interpreted as saying that
the ring $H^*(P_\phi)$ is a (first order) deformation of $H^*(M)$.  
(Here $\phi\in \pi_{d-1}(\Ham(M))$.)   
He  considers only the case when $M= {\CP}^m\times {\CP}^n$ 
and then looks at the part of the vertical quantum product $*_V$ that is
given by counting vertical $A$-curves, where $A$ is the positive
generator of $H_2({\CP}^n,\Z)$.  His equation (5.1) is the analog 
of the statement that  
$$
\io^*: 
QH_V^*(P_\phi, \La_R^*) 
\longrightarrow QH^{*}(M, \La_R^*)
$$
a ring homomorphism with respect to
the quantum products.  The crucial fact in his case is that the moduli
space of {\it spheres} in class $A$ is compact, which means that his operation
gives rise to an associative multiplication on  $\Tilde R\otimes
F[q]/q^2$: cf. property (1$'$) in [Sd2]~Definition 4.1. More generally, it would be
sufficient to require that in the decomposition $$
n_P(u,v,w,z; \io(A)) = \sum_{A = A_1+ A_2} n_P(u, v, e_\al'; \io(A_1))\cdot
n_P(f_\al', w,z; \io(A_2)),
$$
the only nonzero terms have $A_1$ or $ A_2$ equal to
$0$.
\end{rmk}

\subsection{The ring structure of $H^*(P_\phi)$}\label{ssec:ring}

In this section we discuss
conditions under which the ring $H^*(P)$ splits
as the product of $H^*(M)$ with $H^*(S^2)$.  This happens if and only if there is
a ring homomorphism 
$$
r: H^*(M)\to H^*(P),
$$
such that  $\io^*\circ r = {\rm id}$.  Note that in this situation $u_\phi =
r([\om])$ since $r([\om])^{n+1} = r([\om]^{n+1}) = 0$.   Conversely, we have the
following elementary result.

\begin{lemma}  If $H^*(M, \R)$ is generated as a ring by $[\om]$,
then the ring $H^*(P) = H^*(P, \Q)$ splits.
\end{lemma}
\proof{} Since the hypothesis implies that $H^2(M)$ has dimension $1$, we
can assume without loss of generality that $\om$ is rational.  It now suffices
to define the splitting $r$ by $r([\om]^k) = u_\phi^k$.\QED

It is tempting to try to generalize this to the case when $H^*(M)$ is generated
by $H^2(M)$.  However,  in this case the ring $H^*(P)$ does
not always split.  For example, one could calculate $H^*(P)$ as in
Example~\ref{ex:PD} but taking $M$ to be equal to the nontrivial ruled surface
over $S^2$ (rather than $T^2$) so that $H^2(M)$ does generate $H^*(M)$.
See also the examples in [Sd2] which have fiber ${\CP}^m\times {\CP}^n$.

Our next aim is to prove that Theorem~\ref{ring} holds.  Thus we have to
show that the ring $H^*(P)$ splits under the assumption that all vertical
Gromov--Witten invariants 
$$
n_P(v_1,\dots, v_k; B), \;B\ne 0, \; k\le 4,
$$
 vanish. 
Since this implies that the quantum multiplication on $M$ is trivial, we must
have $$
\rho(\phi) = \1\otimes \la + x
$$
as in Lemma~\ref{le:unit}.  Choose $A$ so that the term $\la_A = \mu$ in the
expansion 
$$
\la = \sum_B\la_B e^{-B}
$$
is nonzero.  Then, if $\si_A = \si_\phi + A$,
$$
n_P(a,b, [M]; \si_A) = \mu\, a\cdot b.
$$
The following remark will be useful later.

\begin{lemma}\label{le:inv}  If all the vertical Gromov--Witten invariants
with $B\ne 0$ vanish for $P = P_\phi$, then the same is true for
$P_{\phi^{-1}}$.
\end{lemma}
\proof{}  As noted in the proof of
Corollary~\ref{cor:comp},
 there is a fiberwise diffeomorphism $h$ from $P_\phi$ to
$P_{\phi^{-1}}$  that is a symplectomorphism on the fibers and covers an 
orientation reversing map on the bases.   Hence, given any fibered tame $J$ on
$P_\phi$, there is a fibered tame almost complex structure $J'$ on 
$P_{\phi^{-1}}$  that restricts on the fibers to the pushforward $h_*(J)$.
Therefore, $h$ maps vertical $J$-holomorphic curves on $P$ to vertical
$J'$-holomorphic curves in $P'$.\footnote
{As Seidel pointed out to me, there is no similar correspondence between
curves that represent section classes, since $h$ is antiholomorphic in the
horizontal directions.}
Moreover, as in \S4.3, the regularization $\oMmn =
{\oMm}\,\!^\nu_{0,k}(P,J,B)$ used to calculate the vertical Gromov--Witten
invariants of $P$ can be constructed using this $J$ and a fibered 
perturbation term $\nu$.  Hence the pushforward by $h$ of $\oMmn$
is a regularization of ${\oMm}_{0,k}(P',J',h(B))$ and can be used to calculate
the vertical Gromov--Witten invariants of $P'$.  Thus
$$
n_P(v_1,\dots, v_k; B) = n_{P'}(h_*(v_1), \dots, h_*(v_k); h_*(B)).
$$
 The
result is now immediate. \QED

\begin{lemma} The ring $H^*(P)$ splits if and only if there is an additive
homomorphism
$
s: H_*(M) \to H_{*+2}(P)$ such that, for all $a, b, c\in H_*(M)$,
$$
\begin{array}{lcl}
{\rm (i)} & & s(a)\cap [M] =  a,\\
{\rm (ii)} & & (s(a)\cap s(b))\cdot s(c) = 0.\end{array}
$$
\end{lemma}
\proof{} (i) implies that for each $a,b$ there is $x\in H_*(M)$ such that
$$
s(a)\cap s(b) = s(a\cap b) + \io(x).
$$
Further $s([M]) = [P]$. Hence,
 taking $b = [M]$ in (ii), we find  that
$s(a)\cdot s(c) = 0$ for all $a,c$.   Applying (ii) again, we have
$s(a)\cap s(b) = s(a\cap b)$ for all $a,b$.  One now defines an inclusion $r$ of
the ring $H^*(M)$ into $H^*(P)$
 such that $\io^*\circ r = {\rm id}_M$ by setting
$$
r(\alpha) = \PD_P(s(\PD_M(\al)).
$$
 The converse is similar.\QED

Under the above hypotheses on $P$, define $
s = s_A: H_*(M)\to H_{*+2}(P)
$   by the identity
$$
s_A(a)\cdot_P v = \frac 1{\mu}\, n_P(a, [M], v;\si_A),\qquad v\in H_*(P).
$$
Since
$$
\mu \,(s_A(a)\cap [M])\cdot_M b = \mu\, s_A(a)\cdot_P b =
n_P(a,b, [M];\si_A), $$
it is clear that
$$
s_A(a)\cap [M] = a.
$$
Therefore, we just have to show that
$$
s_A(a)\cdot s_A(b) = 0,\quad (s_A(a)\cap s_A(b))\cdot s_A(c) = 0.
$$
To get a handle on these
intersections, we now show that under the given hypotheses they can be
interpreted in terms of  Gromov--Witten invariants.

\begin{lemma}\label{le:dou} With the  hypotheses of Theorem~\ref{ring},\SS

\NI{\rm (i)} $\mu\, s_A(a)\cdot_P s_A(b) = n_{P,\chi}(s_A(a), s_A(b), [M], [M];
\si_A)$;\SS

\NI{\rm (ii)} $ n_{P,\chi}(s_A(a), s_A(b), [M], [M]; \si_A) = 0$.
\end{lemma}
\proof{} Let $s$ be any splitting that satisfies the conditions of
Lemma~\ref{le:basis}.   Using the composition rule~(\ref{eq:assoc1})  we
find that  
\begin{eqnarray*}
& & n_{P,\chi}(s_A(a), s_A(b), [M], [M]; \si_A)   \\
&  & \qquad\quad =
\sum_{B,i} n_P([M], [M], e_i;\si_A - B)\cdot n_P(s(f_i), s_A(a), s_A(b); B) +\\
& &\qquad \quad \qquad n_P([M],[M], s(e_i);\si_A - B)\cdot 
n_P(f_i,s_A(a), s_A(b); B), \end{eqnarray*}
since all other terms such as 
$$
n_P([M], [M], e_i;B)\cdot n_P(s(f_i), s_A(a), s_A(b); \si_A- B)
$$
that occur in the expansion vanish
by Proposition~\ref{prop:GW}.
By hypothesis, the only nonzero terms have $B = 0$.  For terms of the
first kind,  a dimension count shows that $e_i = [pt]$.  Thus
\begin{eqnarray*}
& & \sum_{i} n_P([M], [M], e_i;\si_A)\cdot n_P(s(f_i), s_A(a), s_A(b); 0)\\
& & \qquad\quad = n_P([M], [M], [pt];\si_A)\cdot n_P([P], s_A(a), s_A(b); 0)\\
& & \qquad\quad = \mu\,
s_A(a)\cdot s_A(b)
\end{eqnarray*}
Further, when $B = 0$ in terms of the second kind, they have to vanish for
reasons of dimension.  This proves (i).

Now let $a = e_j, b = e_k$. We may suppose that  $\dim(e_j)\equiv \dim(e_k)
\equiv \dim(e_i)$$\equiv\eps$ mod $2$ since otherwise all terms below
vanish.  Grouping the terms differently, we have
\begin{eqnarray*}
& & n_{P,\chi}(s_A(e_j), s_A(e_k), [M], [M]; \si_A) \\
&  & \qquad\quad = \sum_{B,i} n_P(s_A(e_j), [M], e_i;\si_A - B)\cdot
n_P(s(f_i),  s_A(e_k), [M]; B) +\\
 & &\qquad \quad\qquad (-1)^\eps n_P(s_A(e_j), [M], s(f_i);B)\cdot n_P(e_i, 
s_A(e_k), [M]; \si_A -B), 
\end{eqnarray*}
where the sign $(-1)^\eps$ is needed because we have reversed the order of
$e_i, s(f_i)$ in the basis. When $B = 0$,
$$
n_P(s(f_i),  s_A(e_k), [M]; B) = n_M(f_i, e_k, [M];B)=(-1)^\eps \,\de_{ik}.
$$ 
Hence  \begin{eqnarray*}
 & & n_{P,\chi}(s_A(e_j), s_A(e_k), [M], [M]; \si_A) \\
& & \qquad\quad =
     (-1)^\eps\,  n_P(s_A(e_j), [M], e_k;\si_A) + (-1)^\eps\, n_P(e_j,
s_A(e_k), [M];\si_A)\\
& & \qquad\quad = \mu\; s_A(e_k)\cdot s_A(e_j) + (-1)^\eps \mu\;
s_A(e_j)\cdot s_A(e_k)\\ & & \qquad\quad = (-1)^\eps 2\mu\, s_A(e_j)\cdot
s_A(e_k), \end{eqnarray*}
where the second equality uses the definition of $s_A$.
This is consistent with (i) only if $s_A(e_j)\cdot s_A(e_k) = 0$ for all $j,k$.
\QED

Hence the splitting $s_A$ satisfies the conditions of Lemma~\ref{le:basis},
and, for simplicity, we now write it as $s(a)$ instead of $s_A(a)$.  

\begin{lemma}\label{le:cap} With the  hypotheses of Theorem~\ref{ring},
$$
n_P(s(a), b, c\,;\si_A) = 0
$$
 for all $ a,b,c\in H_*(M).$
\end{lemma} 
\proof{} By Proposition~\ref{prop:GW}(iv),
\begin{eqnarray*}
 n_P(s(a), b, c\,; \si_A) & = &  \sum_{B,i} n_P(s(a), [M],
e_i; \si_A-\io(B))\cdot n_M(f_i,b,c\,; B) + \\
& & \qquad\qquad n_P([M], [M], e_i; \si_A- \io(B)\cdot n_M(f_i, a, b, c\,; B).
\end{eqnarray*}
But when $B=0$, terms of the first kind vanish by the previous lemma,
since, by definition, $$
n_P(s(a), [M], e_i;\si_A) = \pm s(e_i)\cdot s(a).
$$
Terms of the second kind also vanish  for dimensional reasons,
since $n_P(s(a), b, c\,; \si_A) \ne 0$ only if
$$
\dim(a) + \dim(b) + \dim(c) + 2 = 4n,
$$
which implies that the triple
intersection $(a\cap b)\cdot c = 0$.  Note that here we need to know that
a $4$-point vertical invariant is zero when $B\ne 0$. \QED

\begin{lemma}\label{le:inter} 
With the  hypotheses of Theorem~\ref{ring} and notation as above,\SS

\NI{\rm (i)} $s(a)\cdot_P (s(b)\cap s(c)) = \pm\mu\, n_{P,\chi}(a, s(b), s(c), [M];
\si_A)$;\SS

\NI{\rm (ii)} $ n_{P,\chi}(a, s(b), s(c), [M]; \si_A) = 0$.
\end{lemma}
\proof{} As in Lemma~\ref{le:dou}
\begin{eqnarray*}
 n_{P,\chi}(a, s(b), s(c), [M]; \si_A)   
& = &
\sum_{i} n_P(a, [M], e_i;\si_A)\cdot n_P(s(f_i), s(b), s(c); 0) \\
& &\qquad \quad +  \;n_P(a,[M], s(e_i);\si_A)\cdot n_P(f_i, s(b), s(c); 0).
 \end{eqnarray*}
Terms of the second kind vanish by Lemma~\ref{le:dou}.
Also, because $n_P(a,b,[M];\si_A) $$= \mu a\cdot b$,
 terms of the first kind sum to give
$
(-1)^\eps\mu\; s(a)\cdot_P(s(b)\cap s(c))$ where $\eps = \dim(a)$.  This
proves (i).

We may also write
\begin{eqnarray*}
& & n_{P,\chi}(a, s(b), s(c), [M]; \si_A) \\
&  & \qquad\quad = \sum_{B,i} n_P(a, s(b), e_i;\si_A - B)\cdot n_P(s(f_i),  s(c),
[M]; B) +\\
 & &\qquad \quad\qquad n_P(a, s(b), s(e_i);B)\cdot n_P(f_i,  s(c), [M]; \si_A
-B). \end{eqnarray*}
Again, $n_P(s(f_i),  s(c), [M]; B) = n_M(f_i, c, [M];B)\ne 0$ only if $B=0$.
Since we always have $n_P(a, s(b), e_i;\si_A) = 0$ by
Lemma~\ref{le:cap}, terms of the first kind vanish. A similar argument shows
that terms of the second kind vanish.\QED
 
Here is a more precise version of Theorem~\ref{ring}.

\begin{prop}\label{prop:ring1}   Suppose that
 all vertical $k$-point Gromov--Witten invariants of $P = P_\phi$ with $B\ne
0$ and $k\le 4$ vanish.
 Then:\SS

\NI
{\rm (i)}  the ring $H^*(P)$ splits;\SS

\NI
{\rm (ii)}   $
\oI_c(\phi) =\oI_u(\phi) = 0$;\SS

\NI{\rm (iii)} 
$\rho(\phi) = \mu \1 + x$ for some $\mu\in \Q, x\in QH^+(M)$.
\end{prop} 
\proof{}  The first statement follows from
Lemmas~\ref{le:dou} and~\ref{le:inter}.   Now consider the ring
homomorphism $r: H^*(M)\to H^*(P)$ Poincar\'e dual to $s$.  Thus $$
r(\al) = \PD_P(s(\PD_M(\al)),\quad\mbox{and }\; \io^*\circ r (\al) = \al.
$$
Observe that the coupling class $u_\phi$ must equal $r([\om])$ since 
$r([\om]^{n+1}) = 0$ for dimensional reasons.  Hence 
$$
\PD_P(u_\phi^n) = s(\PD_M(\om^n)) = s([pt]).
$$
Now $s([pt])$ is represented by a union of sections in the (equivalence) class
$\si_A$. Hence it is spherical.  Thus $\PD_P(u_\phi^n)$ is
spherical, so that $\oI_u(\phi) = 0$.  Moreover $u_\phi(\si_A) = 
u_\phi^{n+1}([P]) = 0.$  Since we already know that $c_\phi(\si_A) = 0$ for
dimensional reasons, we must have $\si_A = \si_\phi$. (Note that this
conclusion holds even in the degenerate cases considered in
Lemma~\ref{le:si} (ii).)   Thus (ii) holds and $\rho(\phi)$ has the form
$\mu \1 + x$.\QED

\begin{remark}\rm If there is a fibered $J$ with no vertical
$J$-holomorphic curves at all, one can show in addition that $\mu = \pm 1$.
The point is that the moduli space $\Mm = \Mm(P, J, \si_A)$ of
unparametrized spheres in class $\si_A$ is now a compact manifold and
$\mu$ is the degree of the evaluation map 
$$
\ev: \Mm\times S^2\to P.
$$
($\ev$ is well defined since the elements in $\Mm$ can be parametrized as
sections.)  Hence $\mu $ is an integer.  The proof of Lemma~\ref{le:inv}
shows that if $P = P_\phi$ has no vertical curves, the same is true for $P' =
P_{\phi^{-1}}$.  Since
$\rho(\phi^{-1}) = (1/\mu)\,\1 + x'$, we must have $\mu = \pm 1$.\QED
\end{remark}

The above arguments prove that the ring $H^*(P)$ splits if there are no
nontrivial vertical invariants.  We shall now show that they lead to the same
conclusion either under a strong  hypothesis on the horizontal
invariants or under a ``mixed" hypothesis, that assumes less
about the horizontal invariants but more about the vertical ones.
 Recall from \S 1.3 that, given any section class, the
horizontal quantum product $*_{H,\si}$ is defined by
$$
(u *_{H, \si} v)\cdot w = \sum n_P(u,v,w; \si + B)\otimes e^{-B}.
$$
In particular, $\rho(\phi) = ([M] *_{H,\si_\phi} [M])\cap [M].$  Hence
hypothesis (i) below implies that $\rho(\phi) = \mu \1\otimes e^{-A}$.

\begin{prop}\label{prop:ring2}  The conclusions of
Proposition~\ref{prop:ring1} hold under each of the following assumptions.
\SS

\NI
{\rm (i)}  There is a section class $\si_A$ such that $[M]*_{H,\si_A} [M] =
\mu [P]$; or\SS

\NI
{\rm (ii)} The $3$-point  invariants $n_M(a,b,c\,; B), B\ne 0$ vanish
and $\rho(\phi) = \mu \1\otimes e^{-A}$.

\NI
In particular, in each case  the ring $H^*(P)$ splits as a product of $H^*(M)$
and $H^*(S^2)$. \end{prop}
\proof{} In both cases, one can run through the proofs of
 Lemmas~\ref{le:dou},~\ref{le:cap} and~\ref{le:inter}, checking that the only
nonzero terms in the
expansions have to have $B = 0$.
Then the previous arguments apply.

For example, consider the first expansion in the proof of
Lemma~\ref{le:dou}(i) and suppose that condition (ii) above holds. Then,
because $\rho(\phi) = \mu \1\otimes e^{-A}$,
$$
n_P(a, [M], e_i; \si_A - B) \ne 0 \Longrightarrow B = 0.
$$
Thus $B = 0$ for terms of the first kind.  For terms of the second kind one
uses the identity
$$
n_P(f_i, s_A(a), s_A(b); B) = n_M(f_i, a,b; B)
$$ from Proposition~\ref{prop:GW}(ii). 

Similar arguments work under condition (i) since this is equivalent to
assuming that
$$
n_P([M], [M], v; \si_A - B) \ne 0 \Longrightarrow B=0, v= 
k\,[pt].\qquad\qquad\quad\Box $$

To end this section, we now show that the trivial bundle satisfies
a strong form of condition (i) above that is sufficient to imply
$$
QH^*(M\times S^2) = QH^*(M)\otimes QH^*(S^2),
$$
as claimed in  Proposition~\ref{pr:prod}.  Here
 $QH^*$ denotes the small quantum cohomology ring
described in \S3.2, with coefficients taken either in the usual Novikov ring
$\La_X$ or in the real ring $\La_{R,X}$.   Observe that $\La_{M\times S^2}
\cong \La_M\otimes \La_{S^2}$.

Let $s: H_*(M)\to H_{*+2}(P)$ be the obvious splitting,
$$
s(a) = a\times [S^2].
$$
A proof of the following lemma is sketched in \S4.3.2: see
Lemma~\ref{le:proof}.

\begin{lemma}\label{le:hyp}  Let $P = M\times S^2$ and $\si_0$ be the
section class $[pt\times S^2].$ Then, for all $a\in H_*(M)$ and $u,v,w\in
H_*(P)$,

\NI{\bf
(i)}\,\,
$n_P(u,v,w;D) = 0$ unless $D = k\si_0 + B$ for some $B\in H_2(M)$ and $k =
0,1$.

\NI{\bf (ii)}\,\, $n_P(s(a),u,v;\si_0 + B) = 0$.
\end{lemma}

\begin{prop}\label{prop:prodt} If the conclusions of
Lemma~\ref{le:hyp} hold, then $$
QH^*(M\times S^2) = QH^*(M)\otimes QH^*(S^2).
$$
\end{prop}
\proof{}  By Proposition~\ref{prop:const}, $\rho = \1$ for the
constant loop.  Hence by definition (see Lemma~\ref{le:id}), $n_P(a,b,
[M],\si_O+B) = 0$ if $B\ne 0$ or if $a\cdot b \ne 0$.  It follows easily from
Proposition~\ref{prop:GW} (iv) that
$$
n_P(a,b,c;\si_0 + B) = n_M(a,b,c; B).
$$
Part (ii) of the preceding lemma now implies that the part of
the quantum product in $P$ coming from the
classes $k\si_0 + B$, $k = 0,1$, agrees with that in $QH^*(M)\otimes
QH^*(S^2)$.  But, by (i), classes with $k\ne 0,1$ do not contribute to $QH^*(P)$.
Since they do not contribute to $QH^*(M)\otimes
QH^*(S^2)$ either, the result follows.\QED

\subsection{The nonsqueezing theorem}\label{ssec:nonsq}

First of all, let us calculate the area $\al$ of $(P,\Om)$.
Note that if  $u_S$ is the positive generator of $H^2(S^2,\Z)$ then
$[\Om] = u_\phi + \ka\, p^*(u_S)$, for some $\ka > 0$.

\begin{lemma}  If $[\Om] = u_\phi + \ka\, p^*(u_S)$, then the area of $\al$ of
$(P,\Om)$ equals $\ka$. \end{lemma}
\proof{}  We can write $\Om = \tom + \ka\, p^*(\om_S)$ where $[\tom] =
u_\phi$ and $[\om_S] =
u_S$.  Then $$
\int_P\Om^{n+1} = \int_P (n+1) \ka\, \tom^n \wedge p^*(\om_S) = (n+1)\ka
\int_M \om^n\;\int_{S^2} \om_S,
$$
where the second equality holds because, at each point $x\in P$, $\Om$ restricts to
$p^*(\om_S)$ on the $2$-dimensional subspace $H_x$ of $T_xP$ that is
$\Om$-orthogonal to the tangent space to the fibers.
 The result follows.\QED

\begin{prop}\label{prop:nonsq}  The nonsqueezing theorem holds if 
$$
 n_P([M], [M], [pt];\si_\phi)  \ne 0.
$$
\end{prop}
\proof{}
The hypothesis implies that the evaluation map
$$
\ev_1 : \oMm_{0,2}(P,J,\si_\phi)\to P
$$ 
has nonzero degree, and hence is surjective (where $\ev_1$ evaluates at the
first marked point $z_1$).  Therefore, the usual proof of the nonsqueezing
theorem (see Gromov~[G] or~[LM]) implies that the radius $r$ of any
symplectically embedded ball in $P$ is constrained by the inequality $$
\pi r^2 \le \int_{\si_\phi} \Om.
$$
By the previous lemma, $[\Om] = u_\phi + \al p^*(u_S)$.  Moreover, by
Lemma~\ref{le:si}, we chose $u_\phi$ so that $u_\phi(\si_\phi) = 0$ except in the 
degenerate case when  $\om = 0$ on $H_2^S(M)$  on $H_2^S(M)$.
This shows that $\pi r^2 \le \al$ in all cases except possibly this degenerate
one.

To deal with this case, note that if $\om = 0$ on $H_2^S(M)$ then 
the hypothesis in Proposition~\ref{prop:ring2} holds and
$\si_\phi$  is the Poincar\'e dual of $u_\phi^n$.  Hence we must still have
$u_\phi(\si_\phi) = 0$. \QED

\begin{cor}  The nonsqueezing theorem holds under the hypotheses
of Theorem~\ref{ring} and under those of Proposition~\ref{prop:ring2}.
\end{cor}
\proof{}  By Proposition~\ref{prop:ring1}  $
 n_P([M], [M], [pt];\si_\phi)  \ne 0 $ in these
cases. \QED

\section{Gromov--Witten invariants in fibered spaces}

It remains to prove the properties of the Gromov--Witten 
invariants claimed in Proposition~\ref{prop:GW}
and in Claim~\ref{cl} and Lemma~\ref{le:xx} of \S2.3.  No doubt, any of
the approaches to defining general Gromov--Witten invariants can be
adapted to do this.  We will follow the Liu--Tian~[LiuT1,2] method since in
some way this is the most geometric.  The idea is  to perturb the moduli
space $\oMm_{0,k}(P,J,A)$ to a compact space called
 $\oMmn={\oMm}\,\!_{0,k}^\nu(P,J,A)$ that is sufficiently like a closed
manifold of dimension $d$  to carry a rational fundamental class.  This gives
us a virtual moduli cycle whose elements are stable maps satisfying a
perturbed Cauchy--Riemann equation. We can then argue geometrically using
this cycle.   We have given a detailed description of the construction of
the virtual moduli cycle $\oMmn$ in [Mc3].  In
\S4.1,2 we will try to summarize enough of this discussion to make the 
the rest of the argument in \S4.3 intelligible.

One of the key properties
we need to establish concerns the relation between the invariants in class
$\io(A)$ in  $P$ and those of class $A$ in $M$.  For this, we first need to show
that $\oMmn$ can be constructed using a fibered $J$ and a
perturbation term $\nu$ that is compatible with the fibration $P\to S^2$.
Then  we must  look
carefully at how the stable maps in $M$ sit inside the full cycle 
 $\oMmn$.  In general, they form a codimension $2$ subobject.
We will give a rather {\it ad hoc} definition of the structure put on
 $\oMmn$ that takes this into account: see Definition~\ref{def:s2}.

\subsection{Branched pseudomanifolds}

In this section we describe the topological nature of the virtual cycle.
We will work in
the category of {\it partially smooth} topological spaces $Y$.   Thus
$Y$ is a Hausdorff topological space that is the image of a continous bijection
$$
i_Y: Y_{sm} \to Y,
$$ 
where $Y_{sm}$ is a finite union of open disjoint
(Banach) manifolds.  (Another way to think of $Y$ is as a space
with two topologies.) A morphism $f:Y\to X$ between two objects in this
category is a continuous  map $f:Y\to X$ such that the induced map
$f:Y_{sm}\to X_{sm}$ is smooth.  Usually, we will be interested in
the case when $X$ is a manifold, i.e. when $i_X$ is a homeomorphism, so
that $X$ can be identified with the manifold $X_{sm}$.  Morphisms $f:Y\to
X$ will be called partially smooth maps, or, for short, simply maps.  The path
components of $Y_{sm}$ are sometimes called {\it strata} and sometimes 
{\it components}.  If $Y$ is infinite-dimensional, strata that are open subsets
of $Y$ are called top strata.

\begin{defn}\label{def:brps1}\rm
A (closed, oriented) {\it branched  pseudomanifold}  $Y$ of dimension $d$ 
is a compact partially smooth space $Y_{sm} \to Y$ where the components
of $Y_{sm}$ have dimension at most $d$.  The components of dimension $d$
are called $M_i$, those of dimension $(d-1)$ are called $B_j$, and  $Y_{\le
k}$ denotes the union of all components of $Y_{sm}$ of dimension $\le k$.  
Write:
$$
Y^{top} = \bigcup_i M_i,\quad B = \bigcup_j B_j,\quad 
Y^{sing} = Y - (Y^{top} \cup B) = Y_{\le d-2}.
$$
We assume that each $M_i$ is
oriented and that
$M_i\cup_{j\in J_i} B_j$ 
can be given the structure of a manifold with boundary, where $j \in J_i$ if the
closure $\ov{M}_i$ of $M_i$ in $Y$ meets  $B_j$. 
Similarly, we assume $\ov{B_j} - B_j$ lies in $Y^{sing}$.  Finally, to get a
rational cycle, we assume that $Y$ is {\it labelled}.  This means that the
$M_i$ have rational labels $\la_i\in \Q$, assigned so that the following
condition holds:

\begin{quote}  for each $x\in B$, pick an orientation of $T_xB$ and  divide the
components $M_i$ that have $x$ in their closure into two groups
$I^+, I^-$ according to whether the chosen orientation on $T_xB$  agrees  
with the boundary orientation.  Then we require:
$$
\sum _{i\in I^+} \la_i = \sum_{i\in I^-} \la_i.
$$
\end{quote}
\end{defn}

 Note that in many situations when one is considering branching
 the branching set $B$ is assumed to be closed and so is a codimension $1$
submanifold with singularities.  Here we have considered these singularities
as part of the singular set.  An {\it open} branched pseudomanifold
satisfies all the above conditions except for compactness.  One can also consider
 infinite dimensional
branched pseudomanifolds.  Here the compactness condition is omitted and the
boundary components are  cooriented submanifolds of codimension $1$.

 \begin{lemma}\label{le:cy} Let $Y$ be a (closed, oriented) branched and
labelled 
 pseudomanifold of dimension $d$.
Then every partially smooth map $ev$ from $Y$ to a closed manifold $ X$ 
defines a rational class $ev_*([Y])\in H_d(X)$.  \end{lemma}
\proof{} Since 
we are working rationally, the class $ev_*([Y])$ is defined by
its intersections with classes given by smooth maps $g:Z\to X$ from oriented closed
manifolds $Z$ into $X$.  If $\dim Z + d = \dim X$, then we can jiggle 
$g$ to make it meet each component of the image of $Y_{sm}$ 
transversally.  Hence it will be disjoint from $ev(Y^{sing}\cup B)$ (for
reasons of dimension), and will meet $f(Y^{top})$ transversally in a finite
number of points.  We then set 
$$
 ev\cdot g
= \sum_i \eps_i\la_i, 
$$
where $\eps_i = \pm 1$ is the appropriate sign.
Observe that if $g_0, g_1:Z\to X$ are both transverse to
$ev$ in this way, and are joined by a homotopy  $G:Z\times [0,1]\to X$ then
we can jiggle $G$ (fixing $g_i = G|_{Z\times i}$) so that it is also transverse to
$ev$.  The important point now is that  the image of $G$ still does not meet
$ev(Y^{sing})$ since $Y^{sing}$ has codimension $\ge 2$ though it may now
meet $ev(B)$.    It is not hard to check that the
 condition imposed on the labels in Definition~\ref{def:brps1} guarantees that
this number is independent of the jiggling: cf. [S]. \QED

In order to define the Gromov--Witten invariants  it is sufficient to think of
$\oMmn$ as a (closed, oriented, partially smooth)
 branched pseudomanifold.  However,  in this paper we need a finer
structure in which some strata of codimensions 2 and 3 are controlled.  (Of
course, if necessary, one could control strata of higher codimension.)

\begin{defn}\label{def:s2}\rm  A (closed, oriented, partially smooth)
 branched pseudomanifold $Y$ of dimension $d$ is said to have a
{\it codimension $2$
subcycle}\, $S$ if the following conditions are satisfied:

\NI
(i)   $S$ is (closed, oriented) branched pseudomanifold  of dimension $d-2$;\SS

\NI
(ii)  There is a neighborhood $\Nn(S)$ of $S$ in $Y$ and a commutative
diagram $$
\begin{array}{ccc}
\Nn(S)_{sm} & \to & \Nn(S)\\
\pi'\downarrow & & \pi\downarrow\\
S_{sm} &\to & S
\end{array}
$$
where the maps $\pi'$ and $\pi$ are oriented locally trivial fibrations with
fiber the $2$-disc $D^2$, where $D^2_{sm} = \{0\} \cup (D^2 - \{0\})$.  Here
the components of $\Nn(S)_{sm}$ are given by intersecting $\Nn(S)$ with
the components of $Y_{sm}$.  In particular this means that the top strata in
$S$ are surrounded by top  strata in $Y$ that locally look like  $2$-disc
bundles with the zero section deleted.  A similar remark applies to the
codimension $1$  (branching) components.\SS

\NI
(iii)  It follows from (ii) that the intersection of any top component $M_i$ of
$Y$  with $\Nn(S)$ is the pullback of some top component $M_k^S$ of $S$.  We
require that
$\la(M_i) = \la(M_k^S).$ 
\end{defn}

Here we have made no attempt to put a smooth structure on the 
whole of $\pi^{-1}(S^{top}) = \Nn(S^{top})$ since,
even in the nicest cases, the
disc fibers of $\Nn(S^{top})\to S^{top}$ may have orbifold singularities at their
center: see [Mc2]. The next lemma shows that this is not important.

\begin{prop} \label{prop:codcy} Let $p: P\to S^2$ be a smooth fibration, and
$Y$ be a  branched $d$ dimensional pseudomanifold with codimension $2$
subcycle $S$.  Consider a map $ev:Y\to P$  such that  $S = ev^{-1}(M)$ for some
fiber $M$ of $p$.  Then if $g:Z\to M$ is any cycle such that $\dim Z + d = \dim
P$ 
$$
 ev_*([Y])\cdot_P g_*([Z]) = ev_*([S])\cdot_M
g_*([Z]). 
$$
\end{prop}
\proof{}  We may jiggle $g:Z\to M$ so that it meets $ev(S)$ transversally in $M$
in a finite number of points, and hence define $ ev_*([S])\cdot_M g_*([Z])$.  The left
hand side $ ev_*([Y])\cdot_P g_*([Z])$ is defined by counting intersection points of 
$ev(Y)$ with a jiggling $g'$ of $g$ that meets only the top strata of $Y$.
Since $g$ meets only the top strata of $S$, we are locally in the following situation.
$g$ is a smooth map $Z\to \{0\}\times M\subset \R^2\times M$ and 
$
ev: D^2\times D^{d-2} \to \R^2\times M
$
is a continuous map such that 
\SS

\NI $(a)$  $ev(\{0\}\times D^{d-2}\subset \{0\}\times M;$\SS

\NI
$(b)$  $ev$ is smooth on $\{0\}\times D^{d-2}$ and on $(D^2 - \{0\})\times
D^{d-2}$;\SS

\NI
$(c)$ $ev(\{0\}\times D^{d-2})$ meets $g(Z)$ transversally in $M$ at the single
point  $ev(0,0) = (0, x_0)$.\SS

It follows from $(c)$ that $ev|_{\{0\}\times D^{d-2}}$ and $g$ are local
diffeomorphisms near their intersection point.  Clearly we can perturb $ev$ on
a small neighborhood of $(0,0)$ in 
$D^2\times D^{d-2}$ so that in addition it is smooth in some neighborhood $U$ of
$(0,0)$   and
maps  the $2$-dimensional fibers  in $U$ transversally to $M$.
It is now clear that we can perturb $g$ near $g^{-1}(x_0)$ to a map $g':Z\to P$ 
that 
meets  $ev(D^2\times D^{d-2})$ exactly once transversally (in $P$) and at a point
in the image of $ev((D^2 - \{0\})\times D^{d-2})$.  
The result  then follows since by condition (iii) in Definition~\ref{def:s2} the
labellings of $S$ and $Y$ agree. \QED

\subsubsection{Constructing branched pseudomanifolds} 

In [Mc3] we show in detail how to construct a branched
pseudomanifold from multi-sections of a multi-bundle over a multi-fold. 
The construction amplifies the method of [LiuT1,2].
Here we will give a very brief sketch of the necessary ingredients.

We start with a space $\Ww$ in the partially smooth category that locally
has the structure of an orbifold.  Thus  $\Ww$  is a finite union of
open sets $U_i, i = 1,\dots,k,$, each with uniformizers
$(\TU_i, \Ga_i, \pi_i)$ with the following properties.   
Each $\Ga_i$ is a finite group
acting on $\TU_i$ and the projection
 $\pi_i$ is the composite of the quotient map $\TU_i\to
\TU_i/\Ga_i$ with an identification  $\TU_i/\Ga_i = U_i$.
The inverse image
in $\TU_i$ of each stratum in $U_i$ is an 
open subset of a (complex)  Banach space
 on which $\Ga_i$ acts complex
linearly.   For simplicity, we suppose that the action of $\Ga_i$ is free on
the top strata of $\TU_i$.

The usual compatibility conditions put on the 
local uniformizers of an orbifold are
replaced by the following construction.
 For each subset $I=\{i_1,\dots,i_p\}$ of $\{1,\dots, k\}$  set 
$$
U_I = \cap_{j\in I} U_j,\quad U_\emptyset = \emptyset.
$$  
Let $\Nn$ be the set of all  $I$ for which
$U_I\ne\emptyset$.
For each $I\in \Nn$ define the group
$\Ga_I$ to be the product $\prod_{j\in I} \Ga_j$ and then define $
\TU_I$ to be the fiber product
$$
\TU_I = \{\tx_I=(\tx_j)_{j\in I}:\pi_j(\tx_j) = \pi_\ell(\tx_\ell)\in U_I 
\mbox{ for all
}j,\ell\}\subset \prod_{j\in I} \TU_j,
$$
with its two pullback topologies.
Clearly, $\Ga_I$ acts on $\TU_I$ and the quotient $\TU_I/\Ga_I$ can 
be identified  with $U_I$.  If $J\subset I$, there are projections 
$$
\pi_J^I: \TU_I\to \TU_J, \quad \la_J^I: \Ga_I \to \Ga_J,
$$
where $\pi_J^I$ 
quotients out by the action of the  product group $\Ga_{I-J}$.
We will say that such a collection $$
\TUu = \{(\Tilde U_I, \Ga_I, 
\pi_J^I,\la_J^I):{I\in \Nn}\}
$$
 is a {\it multi-fold atlas} $\TWw$ for $\Ww$ or a multi-fold, for short.
\MS

The next step is to define  a {\it multi-bundle} $\tp:\TEe\to \TWw$.
Suppose given a space $\Ee$
that has two topologies, and a map $p:\Ee\to \Ww$ with the property that each set $E_i =
p^{-1}(U_i)$ has a local uniformizer  $(\TE_i, \Ga_i,\pi_i)$  such that the following
diagram commutes: $$ 
\begin{array}{ccc}
\TE_i & \stackrel{\pi_i}\to & E_i\\
\tp \downarrow & & \downarrow p\\
\TU_i & \stackrel{\pi_i}\to & U_i.
\end{array}
$$
Here we require $\tp: \TE_i\to \TU_i$ to be a $\Ga_i$-equivariant map that
restricts over each stratum of $\TU_i$ to a locally trivial vector bundle.
Thus the fiber $F(\tx_i)$  of $\tp$ at each point $\tx_i$ is a  vector space, but there is
no natural way of identifying one with another  if they lie over points in different strata.
In the application, the points of $\TU_i$ are parametrized stable maps
$\ttau = (\Si, \th)$ and the fiber $F(\ttau)$ is the space
$L^{1,k}(\La^{0,1}(\Si, \th^*TP))$ of sections of a bundle over $\Si$.
When $\ttau$ moves from one stratum to another, the topological type of
$\Si$ may change, and so there is no easy way to identify these fibers.
 Thus, $\Ee$ has
the same local structure as $\Ww$, and we can define a multi-fold atlas $\TEe =
\{\TE_I\}$ for it as above. 

Next we define the concept of a {\it multi-section} $\ts$ of
a multi-bundle $\TEe\to \TWw$.  Very roughly speaking this is a 
compatible collection
$\ts_I$ of (sometimes multi-valued) sections of $\TE_I\to \TU_I$ that in the
case when  $I $ is the singeton $ \{j\}$  is simply an ordinary,
nonequivariant section of the bundle $\TE_j\to \TU_j$. 
The whole point of the construction of $\TEe\to \TWw$ is that
a nonequivariant section $\si(j)$ of $\TE_j\to \TU_j$
can be extended to a multi-section $\ts(j)$ of $\TEe\to \TWw.$  (Here we are
 cheating slightly since one has to refine the multi-fold cover $\{\TU_I\}$ in
order for this to be possible.)  
A multi-section is said to be Fredholm of index $d$ if (among other
conditions) for each $I$ its zero set $\TZ_I$ is a $d$-dimensional
pseudomanifold whose top strata are the intersections with $\TZ_I$ of the
top strata in $\TU_I$.

The second main step in the construction is to
show how to assemble the zero sets of a Fredholm section $\ts$ into a
branched, labelled pseudomanifold $Y_\ts$. There are two steps. First one
replaces the open sets $\TZ_I$ by compact manifolds with boundary
$\TY_I$, and  then one sets
$$
Y_\ts = \coprod_I \TY_I/\sim,
$$
where $\ty_I\sim \tz_J$ if $J\subset I$ and $\pi_J^I(\ty_I) = \tz_J$.
 Elements in the top strata in $Y_\ts$ come
from elements in top strata of $\TU_I$ and
their labels derive from the groups $\Ga_I$.  For example, points in
$Y_\ts^{top}$ that are represented in $\TU_j$ have labels
$1/|\Ga_j|$.  For further details, see [Mc3].

 One of the main questions is how to
construct Fredholm multi-sections.   Below, we will define a
 finite-dimensional family $\ts^\nu, \nu \in R,$ of  multi-sections that are
Fredholm for generic $\nu$.  It follows from the construction that the
homology class represented by the corresponding
branched pseudomanifold $Y^\nu = \oMmn$ is independent of
all choices.

\subsection{The  virtual moduli cycle}

In this subsection we briefly describe how to construct the virtual moduli cycle
$\oMmn = {\oMm}\,\!_{0,k}^\nu(P,J,A)$
as a branched pseudomanifold $Y = Y^\nu$ of dimension $d$, where $d = \dim
P + 2c_1(A) + 2k -6.$  Moreover, when $P$ fibers over $S^2$ and $A$ is a
section class, $Y$ has a codimension $2$ subcycle $S$.  We use the notation of
\S 2.2.1 without further comment and,  for simplicity, will often ignore the
marked points.  In particular, $\Ww$ will denote a suitable neighborhood of
$\oMm$ in the space of all stable maps.

Here are the main features of the structure of the pair $(Y,S)$ when $P$ fibers
over $S$ and $A$ is a section class. The elements of $Y$ are parametrized
stable maps $(\Si, \th)$ that each satisfy a perturbed Cauchy--Riemann
equation  $$
\pJ \th(w)  = \nu_\th (w),\qquad w\in \Si,
$$
where $\nu_\th$ is a $C^\infty$-smooth $1$-form in the Sobolev space
$$
\TL_{\th} = L^{1,p}\left(\La^{0,1}(\Si,
\th^*(T_{vert}P))\right),
$$
consisting of all sections of the  bundle $\La^{0,1}(\Si,\th^*(T_{vert}P))$
that are   $ L^{1,p}$-smooth on each component.  (Here $T_{vert}P$ denotes the
tangent bundle to the fibers of $P\to S^2$.)
There is a finite-to-one map $\pi:Y\to \Ww$ that forgets the
parametrization,  such that each  stratum in $Y_{sm}$ of codimension $2k$
or $2k+1$ is taken to a stratum in $\Ww$ consisting of elements whose
domain has  $\le k+1$ components.  (We would have equality here if it were
not for the extra strata that are introduced to deal with the singular points
of the closure of the branching locus.)  In particular, the top strata $M_i$
and the branch components $B_j$ consist of maps with domain a single
sphere. Further the codimension $2$ subcycle $S$ consists of all elements
whose domain contains more than one sphere.

The evaluation map $ev:\oMmn\to P^k$ is given by
 composing the forgetful map $\pi:Y= \oMmn\to
\Ww$ with the evaluation map $\Ww\to P^k$.  $\oMmn$ is called a virtual moduli cycle
because the class represented by $ev:\oMmn\to P^k$ can be used to calculate
Gromov--Witten invariants.  Presumably it is the same as  the class
$ev_*(\Cc_{0,k}^{hol})$ considered in \S2, but it is not clear whether anyone
has yet checked this.   Therefore from now on we will define Gromov--Witten
invariants using $\oMmn$.  Sometimes we will call $\oMmn$ a {\it regularization} of
$\oMm$.

\subsubsection{The local construction}

The first step is to show that locally $\Ww$  
has the structure assumed in \S4.1.  More
precisely, Liu--Tian show in [LiuT1] Lemma 2.6 that
 each $\tau\in \oMm$ has an open neighborhood $U_\tau$ of the form
$\TU_\tau/\Ga_\tau$,
where  the elements of $\TU_\tau$ are parametrized stable maps 
$\ttau = (\Si,\th)$ and $\Ga_\tau$ acts via reparametrizations.  Thus 
this  local uniformization $(\TU_\tau, \Ga_\tau)$ is given by a constructing a local
$\Ga_\tau$-equivariant slice for the action of the reparametrization group 
on the space of parametrized stable maps.  

Several comments are in order here.  Firstly, one usually wants to consider orbifolds
that are modelled on Banach manifolds.  Thus, best of all would be to have
$\TU_\tau$ an open subset of a Banach space with a linear action of $\Ga_\tau$.
(It is possible, but not easy, to arrange this: see Siebert~[Sb].)  Instead of this, Liu--Tian
show that $\TU_\tau$ is partially smooth.

A second point is that in order to construct their slice Liu--Tian fix the images by $h$
of certain added marked points on $\Si$, which entails the use of the
$  L^{2,p}$-topology on all spaces of maps with domain of fixed topological type. 
This makes gluing more complicated. However, one can avoid this complication (and
instead use the $L^{1,p}$-topology on the strata in $\Ww$)  by Siebert's clever idea of
using integral conditions to fix the parametrization: see [Sb]. With this
approach, certain details of the argument below would have to be changed.

In order to define the overall topology on $\Ww$ one must understand how the
topologies on the  different strata fit together. 
The neighborhood $U_\tau$ is chosen so that all domains $\Si_{\tau'}$ that occur
form a  $2k$-dimensional family (where $\Si_\tau$ has
$k+1$ components)  much as the fibers of the projection $\Ss\to \De$  described
in \S2.3.2.  The parameters in this family are called ``gluing parameters":
their role is briefly described in the discussion just before
Lemma~\ref{le:ztau} below.  Once one knows how to topologize the set of
domains, it is not hard to say what is meant by
nearby stable maps. As a point of notation, we write $\tau = [\Si, h]$ for
elements of $\Ww$ and $\ttau = (\Si, \th)$ for their parametrized lifts.

The next step is to define a local bundle $\TLl_\tau\to \TU_\tau$ whose
fiber at an element $\ttau'= (\Si', \th')$ is the space
$\TL_{\th'}$ defined above.  Note that the fiber changes as when the topological
type of the domain $\Si'$ changes.  But it is a locally trivial bundle over each
stratum and it is not hard to give it a global topology.  Therefore it is a
bundle in our category.  Clearly, the action of  $\Ga_\tau$ lifts to $\TLl_\th$.

 The Cauchy--Riemann operator
$\pJ$ gives rise to a $\Ga_\tau$-equivariant section of $\TLl_\tau \to
\TU_\tau$.  Moreover, its
 linearization at $\ttau'$ is:\footnote
 {If one used Siebert's way of fixing parametrizations one
could consider this as a map from $L^{1,p}$-sections to $L^p$-sections.
Observe also  the domain should be cut down by a finite number of conditions to
fix the parametrizations of the stable maps.  This makes no essential difference to
the present argument, and will be glossed over here.  A precise definition
of the domain of $D\th$ is given in the proof of
Lemma~\ref{le:fibnu} below.}
 $$
D\th':L^{2,p}(\Si', (\th')^*(TP))\to
\TL_{\th'},
$$ 
where $p>2$.
 Now look at the ``center'' point $\tau$ of $U_\tau$ and let $\ttau =
(\Si, \th)\in \TU_\tau$ be one of its lifts. Since $D\th$ is elliptic,
there is  a finite dimensional subspace $R_\tau$ of the vector space
$C^\infty\left(\La^{0,1}(\Si\times P, pr^*TP)\right)$
 such that the operator
\begin{eqnarray}\label{eq:surj0}
D\th \oplus e: L^{2,p}(\Si,\th^*(TP)) \oplus R_\tau \to \TL_{\th}
\end{eqnarray}
is surjective, where $e$ is given by restricting the sections in $R_\tau$
 to the graph of
$\th$.

The next step (see [LiuT1] \S3) is to show that we can choose  $R_\tau$ and a 
partially smooth family of 
embeddings 
$
e_{\ttau'}:  R_\tau  \to \TL_{\th'},
$
so that
\begin{eqnarray}\label{eq:surj}
D\th_{\tau'} \oplus e_{\ttau'}: L^{2,p}(\th_{\tau'}^*(TP)) \oplus R_\tau \to \TL_{\th'}
\end{eqnarray}
is surjective (with uniform estimates for the inverse) for all $\ttau'\in \TU_\tau$. 
Here we are allowed to shrink $U_\tau$. 
 The fact that this is possible is a deep result.  If $\ttau'$ lies in the same
stratum as $\ttau$ it holds because of  the openness of the regularity
condition.    However, to prove this in general 
 one has to show that there is a uniformly bounded family of right inverses to
$D\th_{\tau'}\oplus e_{\ttau'}$ as $\ttau'$ varies over a small enough neighborhood $
\TU_\tau$.  This analytic fact is the basis of all gluing arguments.

With this done, let $pr: \TU_\tau\times R_\tau\to \TU_\tau$ denote the
projection and  consider the pullback bundle $pr^*(\TLl_\tau) \to  \TU_\tau\times
R_\tau.$ This bundle has a section $s$ defined by
\begin{equation}\label{eq:s}
s(\ttau', \nu) = \pJ(\th_{\tau'}) + e_{\ttau'}(\nu).
\end{equation}
By construction, its linearization is surjective  at all points $(\ttau,0)$
and so it remains surjective for $|\nu|\le \eps.$  One then shows by some variant of 
 gluing  that the intersection of $s$ with the zero
section is a (partially smooth) open  pseudomanifold of
dimension equal to ${\rm ind\,}D\th + \dim R_\tau$, whose components are given
by  intersecting the solution space with the different strata in $\TU_\tau$.  By
considering the projection of this solution space to $R_\tau$, one sees that the set of
solutions for fixed generic $\nu\in R_\tau$ is still an open pseudomanifold that we
will call $\TZ_\tau^\nu$.  Since $\TZ_\tau^\nu \subset \TU_\tau$, there is a
forgetful map $\pi:\TZ_\tau^\nu \to \Ww$.

This pseudomanifold $\TZ_\tau^\nu$ has one other very important property
that we will need in order to establish the fact that the virtual moduli cycle
$Y$ has a codimension $2$ subcycle.  Namely, if the domain of $\ttau'\in
\TZ_\tau^\nu$ has $k+1$ components then the nearby elements of
$\TZ_\tau^\nu$ may be obtained from $\ttau'$ by {\it gluing}. When $k=1$
(the case of interest to us) this means the following.  Let $\Ww_2$ denote the
subset of $\Ww$ consisting of elements $\tau$ whose domain has $\ge 2$
components, and suppose that $\TSs$ is a stratum of
$\TZ_\tau^\nu$ that is taken by $\pi$ into a top stratum of
$\Ww_2$.   Then the domain $\Si_{\tau'}$  for $\tau'\in \TSs$   has $2$
components $\Si_1\cup \Si_2$ intersecting at $x= \Si_1\cap \Si_2$.  Hence
there is a line bundle $V\to \TSs$ whose fiber at $\ttau'$ is
$T_x(\Si_1)\otimes T_x(\Si_2)$.  Each sufficiently small element $a$ of the
fiber of $V$ at ${\ttau'}$ serves as a gluing parameter, in the sense that
there is a recipe for making a single sphere $\Si_a$ from the union
$\Si_1\cup \Si_2$ and for constructing a map $g(\ttau', a):\Si_a\to P$.
Moreover   $(\Si_a, g(\ttau',a))\in \TZ_\tau^\nu$  and the map $$
 (\ttau', a)\mapsto (\Si_a, g(\ttau',a))
$$ 
is a bijection from a neighborhood of $\TSs$ in $V$ to a neighborhood of
$\TSs$ in $\TZ_\tau^\nu$. (See [LiuT1] Proposition~3.2 or [FO].)

The above discussion shows:

\begin{lemma}\label{le:ztau} The local solution set $\TZ_\tau^\nu$ has the structure 
of an open branched pseudomanifold with codimension $2$ subcycle
$S = \TZ_\tau^\nu\cap \pi^{-1}(\Ww_2)$.
\end{lemma}

\subsubsection{The global construction}

Because $\oMm$ is compact, it is covered by a finite number $U_j, j = 1,\dots,k$ of the
sets $U_\tau$ and $\Ww$ is defined to be their union.  We write
$(\TU_j, \Ga_j, \pi_j)$ for the corresponding uniformizers, where 
$$
\pi_j:\TU_j\to \TU_j/\Ga_j = U_j
$$
is the obvious projection.  Let $\TLl\to \TWw$ denote the
multi-bundle that restricts on $\TU_j$ to 
$\TLl_j\to \TU_j$.  The Cauchy--Riemann operator $\pJ$, being equivariant,
defines a multi-section of $\TLl\to \TWw$, and the problem  is to define
a suitable family of global perturbations.

To this end, choose a partition of unity $\be_j$ on $\Ww$ subordinate to the
covering $\{U_j\}$  and define the section $\si(j)$ of the pullback  bundle
$pr^*(\TLl_j) \to  \TU_j\times
R_j $ by
$$
\si(j)(\ttau', \nu_j) = \be_j(\tau) \cdot e_{\ttau'}(\nu_j).
$$
Choose $\eps > 0$  so that, for all $j$,  the linearization of $\si(j)$ is
surjective  at all points $|\nu_j|\le \eps$, and put 
$$
R = \oplus_j R_j,\quad R_\eps = \{\nu = (\nu_j)\in R:
 |\nu_j| \le \eps, \mbox{ for all } j.\}
$$
Then set $\Ww_\eps = \Ww\times R_\eps$.  It is shown in [Mc3] that the
pull back multi-bundle $pr^*(\TLl)\to \TWw_\eps$ has a Fredholm
multi-section $\ts$ that is  the sum of $\pJ$ with a finite number
of sections  $\ts(j)$, where $\ts(j)$ restricts over $\TU_j$ to $\be_j\si(j)$.
Hence the corresponding multi-section $\ts\,\!^\nu$ of $\TLl\to
\TWw$ is Fredholm for generic $\nu\in R$.

\begin{prop}\label{prop:cod2}
For generic $\nu\in R$ the zero sets $\TZ_I^\nu$ of the multi-section
$\ts\,\!^\nu$ of  $\TLl\to \TWw$ fit together to give a branched
pseudomanifold $\oMmn$ of dimension $d$ and with codimension $2$
subcycle. \end{prop}
\proof{}  The proof is sketched above except for the statement about the codimension
$2$ subcycle.  Let $\Ww_2$ denote 
the set of all elements in $\Ww$ whose domain has
$> 1$ component, and choose the covering $U_j, j = 1,\dots, k,$  of $\oMm$ so that
its first $\ell$ elements have centers $\tau\in \Ww_2$ and so that the set
$\cup_{j> \ell} U_j$ is disjoint from some neighborhood $\Nn_2$ of
$\oMm\cap \Ww_2$.  It follows from Lemma~\ref{le:ztau} that the local
solution set $\TZ_I^\nu$ of $\ts^\nu$ has a codimension $2$ subcycle 
$$
\TS_I = \TZ_I^\nu\cap \pi^{-1}(\Ww_2)
$$
for $I\subset\{1,\dots,\ell\}$, and for other $I$ does not meet
$\pi^{-1}(\Ww_2)$.  Therefore, we just have to see that in the process of
forming the pseudocycle $Y$ these local subcycles fit together to form a
global subcycle $S$.  The identifications needed to make $Y$ come from the
combinatorics of two subcovers, one a shrinking $\{U_j^0\}$
of $\{U_j\}$ and another associated subcover
$\{V_I\}$ of  $\{U_I\}$: see [LiuT1] Lemma~4.3 and [Mc3] Lemma~4.9.
It follows from the construction detailed in [Mc3]\S4.3 that all we need in
order for the subcycle $S$ to exist is that these subcovers are chosen to be
products near $\Ww_2$.  Thus, if  $r: \Nn_2\to 
\Ww_2$ is  a retraction, it suffices that
$$
V_I = r^{-1}(V_I\cap \Ww_2),\quad U_j^0 =  r^{-1}(U_j^0\cap\Ww_2),
$$
for $j,I\subset \{1,\dots,\ell\}$.
Since it is possible to define the $V_I, U_j^0$ in this way, the result
follows.\QED

\begin{defn}\label{def:allreg}\rm  Any $d$ dimensional branched
and labelled
pseudomanifold $Y = \oMmn$ constructed as above from some covering of 
 a neighborhood $\Ww$ of the space $\oMm_{0,k}(P,J,A)$  will be called a
{\em  regularization} of  $\oMm_{0,k}(P,J,A)$.   \end{defn}

By considering a $1$-parameter version of this construction,
and by considering what happens under a refinement of the covers, 
Liu--Tian show that, for any 
  regularization $\oMmn$  of 
$\oMm_{0,k}(P,J,A)$, the
rational bordism class of the map 
$\ev: \oMmn \to P^k$ is independent of the choice of
$J$, of coverings $\{U_j\}$, and of all other auxiliary structures.  
It depends only on the deformation class of 
of $\Om$ on $P$, on the homology class $A$ 
and on the number of marked points $k$.

\subsection{Properties of Gromov--Witten invariants}

In this subsection we prove Proposition~\ref{prop:GW} and complete the arguments
started in \S 2.4.  This will complete the proof of Theorem~\ref{main}, and hence of
all our other results.  We begin by describing the special properties of
Gromov--Witten invariants in fibered spaces.

\subsubsection{Gromov--Witten invariants in $P\to S^2$}

The following result is immediate from the construction given in \S4.2.

\begin{prop}\label{prop:perturb}  Suppose that $V$ is a subbundle of $TP$ such that
for all representatives $\ttau = (\Si, \th)$ of the elements $\tau \in \oMm$ the cokernel
of $D\th$ is spanned by elements of the space
$$
\Ll_\th^V=L^{1,p}\left(\La^{0,1}(\th^*(V))\right).
$$
 Then we can choose $R$ and the embeddings $e$ so that for all $\nu\in R$ and all
$\ttau$ 
$$
e_{\ttau}(\nu) \in \Ll_\th^V
$$
\end{prop}

We now consider the case when $P = P_\phi$ fibers over $S^2$.
Recall from Definition~\ref{def:Jfib} what it means for $J$ to be compatible with this
fibration.

\begin{lemma}\label{le:fibnu}  Let $P\to S^2$ be a symplectic fibration with compatible
almost complex structure $J$, and suppose that $A\in H_2(P)$ is either a section or
a fiber class. Then we can take $V$ in the above proposition to be the vertical tangent
bundle $T_{vert}P$. \end{lemma}
\proof{}  By Lemma~\ref{le:stembr1}, if $J$ is fibered every component $C_i =
\th(\Si_i)$ of the image of a $J$-holomorphic stable map is either contained in  a
single  fiber,
or is a section. If $C_i$ is a section,
the bundle  $\th_i^*(TP)$ splits as the direct sum $\th_i^*(TC_i)\oplus
\th_i^*(T_{vert}P)$.  Correspondingly 
$$
L^{1,p}\left(\La^{0,1}(\th_i^*(TP))\right) =
L^{1,p}\left(\La^{0,1}(\th_i^*(TC_i))\right)\oplus
L^{1,p}\left(\La^{0,1}(\th_i^*(T_{vert}P))\right). 
$$
Direct calculation shows that $D{\th_i}$ maps
$L^{2,p}(\th_i^*(TC_i))$ onto $L^{1,p}(\La^{0,1}(\th_i^*(TC_i)))$.  
If $C_i$ lies in the fiber, then $\th_i^*(TP)$  splits as 
$\C \oplus \th_i^*(T_{vert}P)$, where $\C$ denotes the trivial bundle.  In this
case also, $D{\th_i}$ maps $L^{2,p}(\C)$ onto  $L^{1,p}(\La^{0,1}(\C))$ because the
domain of $\th_i$ is a sphere. 

Thus the result holds for each component separately.  To deal with $\Si$
as a whole, we must develop some notation to describe its intersection
pattern.  Because the components of $\Si$ are attached according to a tree
graph it is easy to see that we can order them, keeping $\Si_0$ as the stem,
so that the union $\Si_{0}\cup\dots\cup \Si_i$ is connected for each $i$. 
Thus each $\Si_i$ for $i > 0$ is attached to a unique component
$\Si_{j_i}$ for ${j_i} < i$.  Hence there is a point $x_i\in \Si_i$ that is
identified with  $y_{{j_i}}\in \Si_{{j_i}}$.  Then  the
domain of the map $D\th$ is  
$$
L^{2,p}(\th^*(TP)) = \{\oplus_i\, \xi_i \in \oplus_i L^{2,p}(\th_i^*(TP))\;:\;
\xi_{i}(x_{i}) = \xi_{j_i}(y_{j_i})\;\; 1\le i \le \ell\},
$$
and its range is
$$
\oplus_{i=0}^\ell L^{1,p}\left(\La^{0,1}(\th_i^*(TP))\right).
$$
Some remarks are needed here.  Note first that the elements of
$L^{2,p}(\th^*(TP))$ satisfy compatibility conditions  at the double points since
the image of a stable map is always connected.  Secondly, in order to fix the
parametrization of the maps near $\th$ (i.e. in order to construct the local slice
$\TU_\tau$), one does not look at all maps near $\th$ but only at those that satisfy
certain extra normalization conditions.  To do this, Liu--Tian add a minimum number
of marked points $w_i$ to $\Si$ to make each component of the curve $\Si$ stable, i.e.
so that each component has at least $3$ special points.  Then for each $w_i$ they
choose a hypersurface ${\bf H}_i$ in $P$ that meets ${\rm Im \,}\th$ transversally at
$\th(w_i)$ and consider only those maps $\th': \th'(w_i)\in {\bf H}_i$: 
see~[LiuT1] \S2. 
We can choose  ${\bf H}_i$ to be the fiber through $\th(w_i)$ if $w_i$ is on the stem
and  so that the horizontal direction is contained in its tangent space at
$\th(w_i)$ if $w_i$ is on a branch.  The domain of $D\th$ is then cut down by
the corresponding restrictions.   Thus, if there is an added point $w_i$ on the
stem then $\xi_0(w_i)$ must be vertical, and if $w_i$ is on the branch
component $\Si_j$ there is a restriction on the vertical (but {\it not} the
horizontal) part of  $\xi_j(w_i)$.  Let us denote this restricted domain by
${\rm Dom}_{res}(D\th)$.

 In order to prove the lemma we need to see that the ``horizontal part"
of the range 
$$
V_{hor} =   L^{1,p}\left(\La^{0,1}(\th_0^*(TC_0))\right) \oplus \left(\oplus_{i=1}^\ell 
L^{1,p}\left(\La^{0,1}(\C)\right)\right)
$$
is in the image by $D\th$ of the restricted domain
${\rm Dom}_{res}(D\th)$.  It follows from what we
have already proved that, for every element $v=(v_0,\dots, v_\ell)$ of
$V_{hor}$, there is $\xi = (\xi_0,\dots,\xi_\ell)$ such that $v_i =
D{\th_i}(\xi_i)$.  (The reader can check that the extra restrictions put on the
$\xi_i$ by the $w_i$ make no difference here.)
Moreover, we can take $\xi_i\in L^{2,p}(\C)$ for all $i>0$. Thus, when $i > 0$, $\xi_i$
is simply a function on $\Si_i$ and we can alter it by a constant without changing its
image under $D{\th_i}$.  Hence, if we define $\xi' = (\xi_0',\dots,\xi_\ell')$ inductively
by setting 
$$ 
\xi'_0 = \xi_0,\quad \xi_i' = \xi_i + \xi_{j_i}(y_{j_i}) - \xi_i(x_i), $$
we will have
$$
\xi' \in {\rm Dom}_{res}(D\th),\quad D\th(\xi') = Dh(\xi) = v.
$$
Hence result. 
\QED

\begin{defn}\label{def:fibr}\rm  The pair $(J,\nu)$ on $P$ is said to be fibered if $J$ is
compatible with the fibration (as described in Definition~\ref{def:Jfib}) and if
$\nu$ is made as above from sections of the vertical tangent bundle $T_{vert}P$.
\end{defn}

The previous two results imply:

\begin{cor}  If $A$ is a fiber or section class in $P_\phi$ it is 
possible to construct the
regularized moduli cycle $\oMmn = {\oMm}\,\!_{0,2}^\nu(P_\phi,J,A)$ using a fibered
pair $(J,\nu)$.
\end{cor}

\begin{prop}\label{prop:stembr2} If $(J,\nu)$ is fibered and $\si$ is a section class,
then
every element $\tau\in \oMmn = {\oMm}\,\!_{0,2}^\nu(P_\phi,J,\si)$ has a
stem-branch structure.  In other words, one component is a section of
$P_\phi$ and all the others lie in fibers.  The top strata 
in  $\oMmn$ have dimension $d= 2n + 4 + 2c_\phi(\si)$, 
and consist of sections of $P$.  The codimension $2$ strata 
consists of elements with precisely two components, the stem
and a branch consisting of one bubble. 
\end{prop}
\proof{}   Each component $h_i$ of $h$ satisfies an equation of the form
$$
\pJ h_i = \nu_i
$$
where $\nu_i$ is a section of $\La^{0,1}(h_i^*(T_{vert}P_\phi))$.  Hence the map
$\pi\circ h_i: S^2\to S^2$ is holomorphic and the first statement follows as
in Lemma~\ref{le:stembr1}.
The statement about the strata follow firstly from the index formulas
in~(\ref{eq:dim2}) in \S2.2.1, and secondly from the arguments in
Proposition~\ref{prop:cod2}. \QED

\begin{remark}\label{rmk:fam}\rm
Suppose that $X\to B$ is a fibration with fiber $M$, compact base manifold $B$, and
structural group $\Symp(M,\om)$.  Then, if $A\in H_2(X, \Z)$ is a fiber class (i.e. is in the
image of $\io:H_2(M)\to H_2(X)$), 
and if $J$ is an almost complex structure on the vertical tangent bundle $T_{vert}X$
that is tamed by the  symplectic forms on the fibers, the above arguments show
that we can regularize the moduli space $\oMm_{0,k}(X,J,B)$ by a fiberwise
perturbation term $\nu$.  The resulting virtual moduli cycle ${\oMm}\,\!_{0,k}^\nu$
is exactly what is needed to define the fiberwise Gromov--Witten invariants
considered by Le--Ono in~[LO] and Seidel in ~[S2].
Note in particular that, as in Proposition~\ref{prop:stembr2}, each element in 
${\oMm}\,\!_{0,k}^\nu$ has image in a single fiber.
\end{remark}

\subsubsection{The Gromov--Witten invariant  $n_P(\io(a), \io(b);\si)$}

 The Gromov--Witten invariant
$n_P(\io(a),\io(b);\si)$ is the intersection number of the cycle
$$
\ev(\oMm\,\!_{0,2}^\nu(P_\phi,J,\si))
$$
 in $P_\phi\times P_\phi$ with a
generic representative of the class $\io(a\times b)$, where $\io$
is the obvious inclusion and the dimensions are such that transverse
intersections are isolated.  In what follows we will always assume that the
pair $(J,\nu)$ is fibered, and that the dimension condition
$$
\dim (a) + \dim (b) + 2c_\phi(\si) = 2n
$$
is satisfied.  The next lemma shows that we can achieve transversality by   
representing  $\io(a\times b)$ by a cycle of
 the form $\al\times \be$ where $\al, \be$ lie in (distinct) fibers of $P_\phi$.  

\begin{lemma}\label{le:transv} 
If $(J,\nu)$ is fibered, then when evaluating $n_P(\io(a),\io(b);\si)$
 we may assume that the representative of the class $\io(a\times b)$
 has the form $\al\times \be$  where $\al, \be$ are
cycles that represent $\io(a)$ and $\io(b)$ respectively and  
lie in distinct  fibers of
$P_\phi$.  Further all intersection points are transverse and occur with elements
 in the top stratum.  
\end{lemma}
\proof{}  Let $Y = \oMm\,\!_{0,2}^\nu(P_\phi,J,\si)$ where $(J,\nu)$ is fibered.  
Choose generic representatives $\al, \be$ of the classes $\io(a),\io(b)$ that  lie in
distinct generic fibers. As before, since we are working with rational homology we
may assume that $\al,\be$ are the images of smooth closed manifolds.
  We have to show 
that all the intersections of $\al\times \be$ with $\ev(Y^{top})$ are transverse in
$P_\phi\times P_\phi$ and that $\al\times \be$ meets no other elements in the 
image of $\ev$.    This would be obvious for dimensional reasons if it were not for
the fact that we allow only a restricted class of perturbations of $\al, \be$.

First consider the composite
$$
(\pi\times\pi)\circ \ev\circ i_Y:\;\; Y_{sm} \;\to\; Y\;\to\;  P\times P \;\to\; S^2\times
S^2. $$
It is smooth, and so we may choose a regular value  $(\ov{x},\ov{y})\in S^2\times S^2$.
Now let $\al$ (resp. $\be$)  be a generic representative of $\io(a)$ in
$\pi^{-1}(\ov{x})$, (resp. of $\io(b)$ in $\pi^{-1}(\ov{y})$.    Consider an intersection
point $(x,y)$ of $\ev\circ i_Y(Y_{sm})$ with $\al\times \be$, and let $V\subset
T_{(x,y)}P_\phi\times P_\phi$ be the span of the tangent spaces to $\al, \be$ and
$(\ev\circ i_Y)_*(TY_{sm})$.  Since $V$ projects onto  $T_{(\pi(x),\pi(y))} S^2\times
S^2$ by hypothesis, the intersection of $V$ with $T_{vert} P\times T_{vert} P$ has
dimension at most $4n$, with equality only if the dimension of $(\ev\circ
i_Y)_*(TY_{sm})$ is maximal, i.e. if the intersection is with a point in $\ev(Y^{top}-B)$.  
But if the dimension is $< 4n$  we could get rid of this
intersection point by perturbing $\al$ and $\be$ in their respective fibers.  Since
we assumed that $\al, \be$ were generic under such perturbations, this is not possible. 
Hence $(x,y)$ is a transversal intersection point, and lies in $\ev(Y^{top})$ as
claimed.\QED

\NI
{\bf Proof of Proposition~\ref{prop:GW}}

Part (i) follows as  at the end of \S 2.2.2.  (Note that the argument given there is no
longer applicable, since it is not clear whether Siebert's invariants coincide with the
ones we are now considering.)
Part (ii) states that 
$$ 
n_P(a, v, w; B) = n_M(a,v\cap[M], w\cap [M]; B).
$$
Choose representatives $\al, \tbe, \tga$ for $a, v, w$ so that
$\al$ lies in a fiber  $F$ of $P$ that intersects $\tbe$ and $\tga$
transversally.  Let $(J,\nu)$ be fibered, and consider the intersection 
$$
\ev_*({\oMm}\,\!_{0,3}^\nu(P,J,B))\cap (\al\times \tbe\times \tga) \subset P^3.
$$
Since each element of ${\oMm}\,\!_{0,3}^\nu(P,J,B)$ is a stable map with image
in a single fiber of $P$, the only curves that concern us are those lying in $F$. 
Clearly,  by arguing as in Lemma~\ref{le:transv} above, one can
make the intersection transverse by perturbing $\al$ so that it
always remains in some fiber that is transver to $\tbe,\tga$.  Further, it is
clear from the construction of the virtual moduli cycle described in \S 4.1-3
that the subset of  ${\oMm}\,\!_{0,3}^\nu(P,J,B)$ consisting of stable
maps with image in $F$  regularizes  the moduli space of
$B$-curves in $M$ in the sense of Definition~\ref{def:allreg}. Hence there is
a bijective correspondence between the intersection points that contribute to
$n_P(a, v, w; B)$ and those that contribute to
$ n_M(a,v\cap[M], w\cap [M]; B)$.  Since the orientations also correspond,
the result follows.\QED

Part (iii) states that  for all  $v\in H_*(P)$, $v\ne [P]$
$$
 n_P(a, b, v\cap [M], [M]; \si) = n_{P,\chi}(a, b, v, [M]; \si).
$$
To prove this,  construct 
${\oMm}\,\!_4^\nu = {\oMm}\,\!_{0,4}^\nu(P,J,\si)$ as usual, choose $t$ to be a
regular value of the cross ratio
$
cr:\;\;{\oMm}\,\!_4^\nu \to S^2.
$ and set 
$$
{\oMm}\,\!_{4, t}^\nu = cr^{-1}(\chi_0).
$$
As in  Lemma~\ref{le:transv}, one
can make the cycle $\al\times \be\times \nu \times M_0$
transverse to the image of ${\oMm}\,\!_{4,t}$ by perturbing $\al$ and $\be$
within given (generic) fibers $F_0, F_1$ of $P$, by taking $M_0$ to be a
generic fiber $F_\infty$ and by perturbing $\nu$ within the class of
cycles that are products near  the fiber $F_z$  at the point $z$
corresponding to the given cross ratio $t = {\it cr}(0,1,z,\infty)$.
Then it is easy to see that the cycle
$\al\times \be \times (\nu\cap F_z)\times M_0$ is transverse to the image of
${\oMm}\,\!_4^\nu$.  Hence every point counted in
$n_{P,\chi}(a, b, v, [M]; \si)$ is also counted in
$ n_P(a, b, v\cap [M], [M]; \si)$.  Moreover, because the classes $a,b, v\cap [M],
[M]$ are all represented by vertical cycles, the converse also holds.\QED

Part (iv) should be considered as the analog of (iii) for nonvertical classes.
$n_P(v,a,b;\si)$ counts the intersection of $\oMmn = 
{\oMm}\,\!_{0,3}^\nu(P,J,\si)$ with a representative of $v\times a\times b$.
This time, put the representatives $\al, \be$ for $a,b$  in the same fiber. Then
the intersections of $ev(\oMmn)$ with the representing cycle for  $v\times
a\times b$ all occur in the image of the codimension $2$ subcycle $S$ of stable
maps with $\ge 2$ components.   Moreover, we are interested only in that
part $S'$ of $S$ on which the last two marked points $z_2, z_3$ lie on the
same branch.  Then if $(P\times P)_\Delta$ denotes the inverse image of the
diagonal $\Delta$ under the map $P\times P\to S^2\times S^2$, there is an
evaluation map
$$
\ev': S' \to P\times (P\times P)_\Delta,
$$
and it follows from
Proposition~\ref{prop:codcy}, that
$n_P(v,a,b;\si)$ is precisely the number of intersection points 
of $\ev(S')$ with $\nu\times \al\times \be$ when this is considered as a cycle
in $P\times (P\times P)_\Delta$.  There are two kinds of curves counted
here: one with $z_1$ on the stem and the other with $z_1$ on the same branch
as $z_2, z_3$.  The first kind add up to the first sum:
$$
\sum_{B,i} n_P(v,[M], e_i;\si -B)\cdot n_M(f_i, a,b; B)
$$
and the second kind to the second sum:
$$
\sum_{B,i} n_P([M],[M], e_i; \si _ B)\cdot n_M(f_i, v,a,b; B).
$$
The argument here follows the usual proof of the associativity rule for the
quantum product (see~[MS]~\S8.2, for example.)  The point is that one wants
to count intersecting pairs of curves in $P$, one in class $\si -B$ and the other
in class $B$.  Thus one is looking at the intersection of the moduli space of all
such pairs of curves with the diagonal in $P\times P$.  But the diagonal
lies in the homology class $\sum_i e_i\times s(f_i) + s(e_i)\times f_i$.  One
now checks that in both cases the intersection with the classes $s(e_i)\times
f_i$ is zero.  Moreover, the intersection with $e_i\times
s(f_i)$ is in the first case just 
$$
\sum_B n_P(v,[M], e_i;\si -B)\cdot n_M(f_i, a,b; B).
$$
The second case is similar, except that here one needs to use a version
of Proposition~\ref{prop:GW}(ii) that takes account of the fact that 
we are working in $P\times (P\times
P)_\Delta$ instead of $P^3$.  Further details are left to the reader.\QED

\NI
{\bf Proof of Claim~\ref{cl}}

Here is a restatement of this claim in our current language.  The 
statement below makes sense because the map $gr$ extends to a
neighborhood of $\oMm_{0,3}(M,B,J_M)$ in the space of stable maps, and
hence in particular to $\oMm\,\!^\nu_{0,3}(M,B,J_M)$.

\begin{lemma}\label{cl1}  Given any regularization 
$$
\oMmn = 
\oMm\,\!^\nu_{0,3}(M,B,J_M)
$$
 of $ \oMm_{0,3}(M,B,J_M)$,  the invariant
$n_P(a,b;\si_0 + B)$ is the intersection number of $ev\circ
gr(\oMmn)$ with a generic cycle in $M\times M$ representing
$a\times b$.
\end{lemma}
\proof{} Fix two fibers $M_0, M_1$ of $P$ and let $\Ww'$ be the space   of
stable maps $\tau = [\Si, h, z_0, z_1]$ for which $h(z_i)\in M_i$.  We claim  that
the branched pseudomanifold $gr(\oMmn)$   regularizes 
the moduli space
$$
\Mm'(P) = \oMm_{0,2}(M\times S^2, \si_0+\io(B)\cap \Ww'.
J)
$$ 
To see this, first observe that the expected or formal dimension of $\Mm'(P)$ is
$$
2n+2 + 2c_1(\si_0+B) - 6 = 2n + 2c_1(B) = \mbox{f. dim.}\,
\Mm_{0,3}(M,B,J), 
$$
where  $6 = \dim \PSL(2,\C)$ is subtracted because, if the marked points
lie on the stem component and if this is parametrized as a graph, the position
of the marked points is uniquely determined.  (As a check, note also that  if a
marked point $z_i$ is in a branch component, this component is constrained to
lie in
 the chosen fiber $M_i$  and so the two extra dimensions from
the possible movement of $z_i$ are cancelled out by this codimension $2$
constraint on the branch.)  Hence $gr(\oMmn)$ does have the right dimension. 
Next, observe that the images by $gr$ of the choices (for example of coverings
$U_i, \TV_I$) made in the construction of $\oMmn$  are valid choices  for the
construction of a regularization of $\Mm'(P)$.   This proves the claim.

It remains to check that
$n_P(\io(a), \io(b);\si)$ can be calculated from this
regularization via the  evaluation map into $M_0\times M_1$.
But this was proved in Lemma~\ref{le:transv}.\QED

This lemma immediately implies that $n_P(a,b;\si_0 + B)$  is always
$0$.  The basic reason is that $ev\circ gr_*(\oMmn)$ reduces dimension:
$\dim(ev\circ gr_*(\oMmn))$ is strictly smaller than $\dim(\oMmn) = 2n +
2c_1(B)$ since the image of $\tau$ does not change when the marked point 
$z_2$ is moved in its given component.  
To be more precise, note that  $n_P(a,b;\si_0 +
B)\ne 0$ only when $$
\dim(a) + \dim(b) + 2c_1(B) = 2n.
$$
Further the intersection points of $ev\circ gr_*(\oMmn)$ with $a\times b$
are precisely those of 
$$
ev_M(\oMmn\,\!_{0,2}) \to M\times M
$$
 with $a\times b\times M$.
But the image of $ev_M$ has codimension $2n - 2c_1(B)- 2$ in $M\times M$
and so  one can move $a\times b$
 to be disjoint from this image.\QED

We end this subsection by sketching a proof of Lemma~\ref{le:hyp}
which says:

\begin{lemma}\label{le:proof}  Let $P = M\times S^2$ and $\si_0$ be the
section class $[pt\times S^2].$ Then, for all $a\in H_*(M)$ and $u,v,w\in
H_*(P)$,

\NI{\bf
(i)}\,\,
$n_P(u,v,w;D) = 0$ unless $D = k\si_0 + B$ for some $B\in H_2(M)$ and $k =
0,1$.

\NI{\bf (ii)}\,\, $n_P(s(a),u,v;\si_0 + B) = 0$.
\end{lemma}
\proof{}  The proof of (ii) is very similar to that of Lemma~\ref{cl1}.
Let us first consider 
$$
n_P(a, s(b), s(c);\si_0 + B).
$$
This time, because only one of the classes $a, s(b), s(c)$ is in $H_*(M)$,
we fix one fiber $M_0$ and let  $$
\Ww' = \{\tau= [\Si, h, z_0, z_1,z_2]\;:\; h(z_0)\in M_0\}.
$$
Further define
$$
gr': \oMm_{0,5}(M,B,J_M) \to \oMm_{0,3}(M\times S^2, \si_0 + B, J)\cap
\Ww' $$
by setting
$$
gr'([\Si, h, z_0, \dots, z_4]) = [\Si, \widehat{h}, z_0, z_3, z_4]$$
where,
as in the definition of $gr$,
$$
\widehat{h} (z) = (h(z), \ga(z))\in M\times S^2
$$
 is defined using the
first three marked points $z_0,z_1,z_2$.  As before, it is not hard to check
that if $\oMmn$ regularizes  $\oMm_{0,5}(M,B,J_M)$, then $gr'(\oMmn)
$ regularises $\oMm_{0,3}(M\times S^2, \si_0 + B, J)\cap
\Ww' $ in $\Ww'$.  Moreover, there is an evaluation map
$$
ev'\;:\: gr'(\oMmn) \to M\times P \times P.
$$
Clearly, $ev'(gr'(\tau))$ intersects  $ a\times s(b)\times s(c)$
exactly when $h(z_0)\in a, h(z_3)\in b, h(z_4)\in c$.
But $n_P(a, s(b), s(c);\si_0 + B) \ne 0$ only when
$$
2n+2 + 4 + 2c_1(B) + \dim(a) + \dim(b) + \dim (c) + 4 = 6n+6.
$$
i.e. if $2c_1(B) + \dim(a) + \dim(b) + \dim (c) = 4n - 4$.
But in these dimensions there are for generic $J$ no $J$-curves in $M$
through $a,b,c$.  More precisely, representatives for $a,b,c$ 
(still called $a,b,c$) can be chosen 
in $M$ so that the image of
$$
ev:\oMm_{0,3}(M,B,J_M)\to M\times M\times M
$$
does not meet $a\times b\times c$.  

Similar arguments show that $n_P(s(a),b,c;\si) = n_P(s(a), s(b), s(c);\si) = 0$
whenever $\si = \si_0 + B$.

We now have to see that $n_P(u,v,w; k\si_0 + B) = 0$ whenever $k\ne 0,1$.
Note first that these invariants vanish when $k < 0$, since positivity of 
intersections with the fibers $[M]$ of $P$ implies that  the moduli spaces
$\oMm_{0,3}(P, k\si_0+B, J_M\times j)$ are empty when $k < 0$.
When $k>1$ the idea again  is to show that the evaluation map
reduces dimension by constructing elements of
$\oMm_{0,3}(P, k\si_0+B, J_M\times j)$ from elements of $\oMm(M, B,
J_M)$ and $\oMm(S^2, k \si_0, j)$. (By slight abuse of
notation we have written
 $\si_0$ for the fundamental class of $S^2$.)
   The
only problem is to find a replacement for the map $gr$.

Consider the map
$$
\pi = \pi_M\times \pi_S:  \oMm_{0,3}(M\times S^2, k\si_0 + B, J_M\times
j) \to \oMm_{0,3}(M, B, J_M)\times 
\oMm_{0,3}(S^2, k\si_0, j)
$$
obtained by projection followed by contraction of unstable components,
and let $\Nn$ be its image.  This map $\pi$ is injective on the top
component,
but need not be injective in general.\footnote
{To get injectivity, one must also keep
track of information about the domain $\Si$ of the element in the left hand
side and how it relates to the domains of its projections to $M$ and $S^2$.
For example, suppose that $\Si$ has three components $\Si_1$, that maps
into a fiber, and $\Si_i, i = 2,3$ where $\Si_i$ maps to the flat section
$p_i\times S^2$.  Then the image of this element $\pi$ is independent of the
choice of $p_2,p_3\in M$.}

It is easy to check that for each stratum $\Ss$ in
$\Nn$ the restriction 
$$
\pi: \pi^{-1}(\Ss) \to \Ss
$$
is a fibration.  Note also that  $\oMm_{0,3}(S^2, k\si_0, j)$ is an orbifold of
the correct dimension and so is already regular.  (All the relevant maps
$D\th$ are already surjective.) 
Consider the pullback 
$$
\oMmn(M\times S^2) = \pi^{-1}(\oMmn(M)\times \oMm_{0,3}(S^2, k\si_0, j)
$$
where $\oMmn(M)$ 
is a regularization  of $\oMm_{0,3}(M, B, J_M)$.  By the above remarks this is
a branched pseudomanifold $Y$.  One now has to check that $Y$ is homologous
to a regularization of $\oMm_{0,3}(M\times S^2, k\si_0 +B, J_M + j)$.
Note that $Y$ is made in much the same way as is the regularization.  However
it is made from a cover of the moduli space that is pulled back from the
product of a cover of  $\oMm_{0,3}(S^2, k\si_0, j)$ with $\oMm_{0,3}(S^2,
k\si_0, j)$, and these pullback sets are larger than the sets used in the
constructions in \S 4.2.   However, these constructions  behave well
when the cover is refined, and, using this, the desired result follows fairly
easily.  

Granted this, we can use this pullback set $\oMmn(M\times S^2)$
when evaluating the Gromov--Witten invariants for $M\times S^2$.
Now observe that when $k > 1$ the image of the evaluation map
$$
\oMm_{0,3}(S^2, k\si_0, j)\to (S^2)^3
$$ 
has dimension strictly smaller than 
$$
\dim(\oMm_{0,3}(S^2, k\si_0, j)) = 2 + 4k.
$$
Hence, for dimensional reasons, this evaluation map never gives rise to
nonzero Gromov--Witten invariants.  This implies that the evaluation map
$$
\oMmn(M)\times 
\oMm_{0,3}(S^2, k\si_0, j)\;\stackrel{ev_M\times ev_S}\longrightarrow\;
M^2\times (S^2)^3 $$
also has too small an image to give rise to nonzero 
 Gromov--Witten invariants.  Hence the same is true for 
$$
ev_{M\times S}:\oMmn\,\!_{0,3}(M\times S^2)\to (M\times S^2)^3
$$
since this map factors through $ev_M\times ev_S$.\QED

\subsubsection{The composition rule}

It remains to prove Lemma~\ref{le:xx}  which is the last step in the
proof of Proposition~\ref{prop:comp}.  We assume that we are in the
situation described in \S 2.3.2. The space $\Xx$ fibers over the disc $\De$
with fibers $X_t$ that are themselves fibered.  Thus the projection   $\Xx\to
\De$ factors as $$
\Xx\stackrel{q}\longrightarrow \Ss\stackrel{\pi}\longrightarrow \De,
$$
and the construction of the perturbation term $\nu$ must now be compatible with 
both fibrations.   Recall that $\Ss$ is a holomorphic space and that $\pi:\Ss\to \De$ is
  holomorphic with a singular fiber at $0$.  This singular
fiber is the transverse union of two $-1$ spheres.  Hence all these fibers in $\Ss$ are
regular  as far as Fredholm theory goes.  It follows easily that the proof of
Lemma~\ref{le:fibnu} can be adapted to show:

\begin{lemma}\label{le:fibnu2}  In this situation, we can take the bundle $V$ in
Proposition~\ref{prop:perturb} to be  the vertical tangent bundle
$T_{vert}\Xx$ of the fibration $\Xx\to \Ss$.
\end{lemma}

Note that the fiber of $T_{vert}\Xx$ at the point $x\in X_t$ can be identified with the
fiber of $T_{vert}(X_t)$ at $x$, in other words it consists of vectors tangent to $M$.

We now construct local uniformizers $(\TU_\tau,\Ga_\tau)$ 
and families of embeddings $e_{\ttau'}: R_\tau\to \TL_{\th'}^V$ for $\ttau'\in
\TU_\tau$ as usual, starting first with the points $\tau\in \Ww(X_0)\subset
\Ww(\Xx)$.   The aim is to do this compatibly with the projection $\Xx\to \De$,
i.e. so that the regularization $Y(\Xx)$ can be thought of 
as a family of regularizations
$Y(X_t)$, for $t$ near $0$.
To do this one just follows the construction given in 
Proposition~\ref{prop:cod2} and Lemma~\ref{le:ztau} for a codimension $2$
subcycle. In other words, make a regularization $Y(X_0)$ and then extend
this over a neighborhood in $\Ww(\Xx)$ by gluing.  The family of curves
$\Si_t$ is precisely what one obtains by gluing the singular $\Si_0$ with
parameter $t$.  Hence, because the perturbation term $\nu$ projects to $0$
in $\Ss$,  when one  glues a lift $C_0\subset X_0$ of $\Si_0$   with
parameter $t$  one gets a  curve $C_t$  in $X_t$.

Therefore, one can construct $Y(\Xx)$ so that,
for all
sufficiently small $t$,
 the inverse image
$Y(X_t)$ of $\Ww(X_t)$ in $Y(\Xx)$ is a regularization for $\oMm(X_t)$.  Of course
$Y(\Xx)$ is not compact, but the restrictions $Y(X_t)$ are all compact. \MS

\NI{\bf Proof of Lemma~\ref{le:xx}}

We have now shown how to construct a regularization  $Y(\Xx)$ of 
$\oMm_{0,4}(\Xx,\Tilde J, \Tilde \si)$ that is compatible with the fibration $\Xx\to
\De$.  Let $Y_{\chi}(\Xx)$ be the  subset  on which the cross ratio is fixed.
As before $Y_\chi(\Xx)$ is a union of virtual moduli cycles $Y_\chi(X_t)$.

The first step is to prove that one can achieve transverse 
intersections of $\ev (Y_\chi(\Xx))$ with the class
$\prod ([\Tilde \al_i])$ by choosing generic cycles $\Tilde\al_i$ that are compatible
with the fibrations on $\Xx$ as described in \S2.3.2.  The needed arguments are very
similar to those in Lemma~\ref{le:transv} and will be left to the reader.

Next observe that because of the special form of the cycles $\Tilde \al_i$  the only
elements in $Y_{\chi}(\Xx)$ that meet the cycles $\Tilde\al_i$ lie in a single set
$Y_\chi(X_t)$.
 Hence one gets an intersection of
$\ev(Y_{\chi}(X_t))$ with the cycle $\prod(\Tilde\al_i)_t\in (X_t)^4$, where $(\Tilde
\al)_t = \Tilde \al\cap
X_t$.  However this intersection is not transverse, since it can be destroyed by
moving  the fiber $(\Tilde\al_2)_t = M$ to another fiber in $X_t$.  If one did this in
$\Xx$ the effect would be to change the relevant value of $t$.  
It is now not hard to see that one can get  a transverse
intersection in $X_t$ by considering the intersection of $\prod(\Tilde\al_i)_t$ with the
image of the whole of $Y(X_t)$ rather than just $Y_\chi(X_t)$.  
This proves (i).

Now consider (ii).   Note that the cross-ratio of the $4$ marked points is fixed on the
component of $Y(X_0)$ on which  there are two marked points mapping to each
component of $\Si_0$. Since  this is the only component that intersects the cycle
$\prod(\Tilde\al_i)_0$ , we can think that $Y_\chi(X_0) = Y(X_0)$.
Moreover, again we can achieve transversality of intersection in $X_t$ for all $t$ near
$0$ by taking generic cycles $\Tilde\al_i$ of the given form. 

Because $Y(X_0)$ is a codimension $2$ subcycle
in $Y(\Xx)$, it follows from Proposition~\ref{prop:codcy} that the intersection number
$n_{X_0}(a,[M],[M],b)$ of $Y(X_0)$ with $\prod(\Tilde\al_i)_0$ is the same as the
intersection number $n_{X_t}(a,[M],[M],b)$ of $Y(X_t)$ with
$\prod(\Tilde\al_i)_t$  for each fixed $t$ near $0$.  It remains to check that 
$n_{X_0}(a,[M],[M],b)$ is
precisely $$
 \sum_{\si = \si_1\#\si_2 }n_{P_\phi}(a,[M],e_i;\si_1)\;
n_{P_\psi}(f_i,[M],b;\si_2).
$$
This completes the proof of Lemma~\ref{le:xx}.  Note that a very similar argument
would prove the composition rule~(\ref{eq:assoc1}) for general
Gromov--Witten invariants.\QED

 %%%%%%%%%%%%
\MS

 \NI
{\large\bf References}
\MS\MS

\NI
[An]  S. Anjos,  in preparation.\MS

\NI
[B]  A. Blanchard,  Sur les vari\'et\'es analytiques complexes, {\it Ann. Sci.
Ec. Norm. Sup. }, {\bf (3) 73},  (1956).\MS

\NI
[FO]  K. Fukaya and K. Ono, Arnold conjecture and
Gromov--Witten invariants, preprint (1996)
\MS

\NI
[Gi] A. Givental, Equivariant Gromov--Witten invariants, {\it Internat. Math.
Research Notes}, (1996), 613--663.

\MS

\NI
[G]  M. Gromov,  Pseudo holomorphic curves in symplectic manifolds, {\it
Inventiones Mathematicae\/}, {\bf 82} (1985), 307--47.

\NI
[HLS] H.~Hofer, V.~Lizan and J.-C.~Sikorav,   On genericity for complex curves in
$4$-dimensional almost complex manifolds, {\it Duke Math Journal},

\MS

\NI
[HS] 
H. Hofer and D.A. Salamon,  Floer homology and Novikov rings,  
{\it Floer Memorial volume\/}  eds. Hofer,
Taubes, Weinstein, Zehnder, Birkh\"auser, Basel, (1996).
\MS

\NI
[HS] 
H. Hofer and D.A. Salamon,  Marked Riemann surfaces of genus $0$,
preprint (1997). \MS

\NI
[IP] E. Ionel and T. Parker, Gromov--Witten invariants of symplectic sums, (in
preparation) \MS

\NI
M. Kontsevich and Yu. Manin,  Quantum Cohomology of a Product.
\MS

\NI
[LMP1] F. Lalonde, D. McDuff and L. Polterovich, On the Flux
conjectures, 
CRM Proceedings and Lecture Notes vol 15, (1998), 69--85

\MS

\NI
[LMP2]  F. Lalonde, D. McDuff and L. Polterovich,
Topological rigidity of Hamiltonian
loops and quantum homology, to appear in {\it Invent. Math}.
\MS

\NI
[LMP3]  F. Lalonde, D. McDuff and L. Polterovich,
 On Hamiltonian fibrations, in preparation.
\MS

\NI
[LO] Hong-Van Le and Kaoru Ono, Topology of symplectomorphism groups and
Pseudoholomorphic curves, preprint (1997)\MS

\NI
[LiT]  Jun Li and Gang Tian, Virtual moduli cycles and 
Gromov--Witten invariants for general symplectic
manifolds, preprint alg-geom 9608032.
\MS

\NI
[LiuT1] Gang Liu and Gang Tian, Floer homology and Arnold
conjecture, {\it Journ. Diff. Geom} (1998).
\MS

\NI
[LiuT2] Gang Liu and Gang Tian, On the equivalence of multiplicative
structures  in Floer Homology and Quantum Homology, preprint 
Oct.1998. \MS

\NI
[Mc1] D. McDuff, From deformation to isotopy, to appear in Proceedings of the 1996
Conference at Irvine, ed. Stern, International Press.\MS

\NI
[Mc2] D. McDuff, Almost complex structures on $S^2\times S^2$, 
 SG/9808008, to appear in {\it Duke Math. Journal}.\MS

\NI
[Mc3] D. McDuff, The virtual moduli cycle, preprint (1998).
\MS

\NI
[Mc4] D. McDuff, Fibrations in Symplectic Topology,  {\it
Documenta Mathematica}, Extra Volume ICM 1998, Vol I, 339--357.
\MS

\NI
[MS1]   D. McDuff and D.A. Salamon, {\it $J$-holomorphic curves and
quantum cohomology}, Amer Math Soc Lecture Notes \#6, Amer.
Math. Soc. Providence (1994).
\MS

\NI
[MS2]   D. McDuff and D.A. Salamon,  {\it Introduction to Symplectic Topology}, 2nd
edition, OUP, (1998).\MS

\NI
[PSS] S. Piunikhin, D. Salamon and M. Schwarz, Symplectic Floer--Donaldson
theory and Quantum Cohomology, {\it Contact and Symplectic Geometry} ed C.
Thomas, Proceedings of the 1994 Newton Institute Conference, CUP, Cambridge (1996)
\MS

\NI
[P] L. Polterovich: Symplectic aspects of the first eigenvalue, preprint (1997)
\MS

\NI
[R] Y. Ruan,  
Virtual neighborhoods and pseudoholomorphic curves, preprint
1996
\MS

\NI
[RT]
Y. Ruan and G. Tian,  A mathematical theory of quantum cohomology.
 Journ Diff Geo {\bf 42} (1995), 259--367.
\MS

\NI
[Sd1] P. Seidel, $\pi_1$ of symplectic
automorphism groups and invertibles in quantum cohomology rings, 
{\it Geometric and Functional Analysis}, {\bf 7} (1997), 1046--95.
\MS

\NI
[Sd2]  P. Seidel,  On the group of symplectic automorphisms of ${\CP}^m\times
{\CP}^n$, preprint (1998) \MS

\NI
[Sb] B. Siebert, Gromov--Witten invariants for general symplectic
manifolds, preprint (1997)  
\MS

\end{document}